\documentclass[12pt]{article}
\usepackage{epsfig, color}
\usepackage{amssymb,latexsym}
\usepackage{hyperref}
\usepackage{epstopdf}
\usepackage{amsfonts,amsmath, mathrsfs}
\usepackage{enumerate}

\usepackage{amsmath} 
\allowdisplaybreaks

\newtheorem{lem}{Lemma}[section]
\newtheorem{theo}{Theorem}[section]
\newtheorem{pro}{Proposition}[section]
\newtheorem{prob}{Problem}[section]
\newtheorem{cor}{Corollary}
[section]
\newtheorem{con}[prob]{Conjecture}
\newtheorem{que}{Question}[section]
\newtheorem{exa}{Example}[section]

\newtheorem{problem}
{Problem}[section]

\newcommand {\red} {\textcolor{red}}

\newcommand {\blue} {\textcolor{blue}}

\usepackage{tikz,chemfig}

\usetikzlibrary{arrows.meta,decorations.markings}

\usetikzlibrary{calc}

\usetikzlibrary{arrows.meta}

\usetikzlibrary{backgrounds}
\usetikzlibrary{shadings,intersections}
\usetikzlibrary{arrows}
\usetikzlibrary{shapes,shapes.geometric,shapes.misc}
\pgfdeclarelayer{edgelayer}
\pgfdeclarelayer{nodelayer}
\usetikzlibrary{automata}
\pgfsetlayers{background,edgelayer,nodelayer,main}
\tikzstyle{none}=[inner sep=0mm]
\tikzstyle{every loop}=[]

\tikzstyle{dottedge}=[dash pattern=on \pgflinewidth off 2pt]
\tikzstyle{dashed}=[dash pattern=on 3pt off 3pt]

\tikzstyle{new style 0}=[fill=black, draw=black, shape=circle]
\tikzstyle{red style 1}=[fill=red, draw=black, shape=circle]
\tikzstyle{blue style 2}=[fill=blue, draw=black, shape=circle]
\tikzstyle{white style 4}=[fill=white, draw=black, shape=circle]
\tikzstyle{bklack style 5}=[fill=black, draw=black, shape=rectangle]
\tikzstyle{red style 3}=[fill=red, draw=black, shape=rectangle]
\tikzstyle{yellow style 7}=[fill=yellow, draw=black, shape=rectangle]
\tikzstyle{new style 8}=[fill={rgb,255: red,0; green,132; blue,0}, draw={rgb,255: red,0; green,131; blue,0}, shape=circle]

\tikzstyle{new edge style 0}=[-]
\tikzstyle{new edge style 1}=[-, draw=red]
\tikzstyle{new edge style 2}=[-, draw=blue]
\tikzstyle{new edge style 3}=[-, draw={rgb,255: red,0; green,156; blue,0}]

\def \seta {{\cal A}}
\def \EC {{\cal EC}}

\newcommand\Enum[1]
{
	\begin{enumerate}
		#1
	\end{enumerate}
}

\newcommand\Eqnn[1]
{
	\begin{eqnarray*}
		#1
	\end{eqnarray*}
}

\newcommand\Eqn[2]
{
	\begin{eqnarray}\label{#1}
		#2
	\end{eqnarray}
}

\newcommand\equ[2]
{
\begin{equation}\label{#1}
#2
\end{equation}
}

\newcommand\Bra[1]
{
	\left ( 	#1	\right )
}

\newcommand{\solu}{{\noindent {\em Solution}.\quad}}

\newcommand{\proof}{{\noindent {\em Proof}.\quad}}
\newcommand{\proofend}{{\hfill$\Box$}\vspace{0.3 cm}}

\usepackage{bm}    

\newcommand{\gap}
{\vspace{0.4 cm}}

\newcommand{\sgap}{\vspace{0.3 cm}}

\def\spose#1{\hbox to 0pt{#1\hss}}
\def\ltapprox{\mathrel{\spose{\lower 3pt\hbox{$\mathchar"218$}}
 \raise 2.0pt\hbox{$\mathchar"13C$}}}
\def\gtapprox{\mathrel{\spose{\lower 3pt\hbox{$\mathchar"218$}}
 \raise 2.0pt\hbox{$\mathchar"13E$}}}

\newcommand{\Z}{{\mathbb Z}}

\newcommand{\R}{{\mathbb R}}

\def \setb {{\cal B}}
\def \W {{\cal W}}
\def \F {{\cal F}}
\def \O {{\cal O}}
\def \seti {{\cal I}}
\def \setg {{\cal G}}
\def \setv {{\cal V}}
\def \sete {{\cal E}}
\def \P {{\cal P}}
\def \G {{\cal G}}
\def \D {{\cal D}}
\def \setc {{\cal C}}

\def \tsim {\stackrel T\sim}

\def \Sorder {\overline {\Omega}}
\def \order {{\Omega}}

\def\tchi {\tilde \chi}
\def \sets {{\cal S}}
\def \F {{\cal F}}
\def \angx {\langle x\rangle}
\def \P {{\cal P}}
\def \OP {{\cal OP}}

\newcommand\stirone[2]
{\left [ #1  \atop #2 \right ]
}

\newcommand {\relabel}[1] {\label{#1}[\red{**its label}: \blue{\bf #1}]}\newcommand {\rebibitem}[1] {\bibitem{#1}[\red{**its label}: \blue{\bf #1}]} 

\def\relabel {\label} \def\rebibitem {\bibitem} 

\renewcommand \red {}
\renewcommand \blue {}

\begin{document}
	
	\baselineskip 0.6 cm

\begin{title}
{A Survey of Graph Polynomials Related to the Chromatic Polynomial}
\end{title}

\author{Fengming Dong\thanks{Emails: fengming.dong@nie.edu.sg and donggraph@163.com.}\\
National Institute of Education\\
Nanyang Technological University}
\date{}
\maketitle{}

\begin{abstract}
	
	Introduced by Birkhoff in 1912 to study the four-color conjecture, the chromatic polynomial has become one of the most extensively studied graph polynomials. It is closely connected to the Tutte polynomial through specialization, to the characteristic polynomial of the associated graphic matroid, and to the flow polynomial through planar duality. It has also inspired several related polynomials, including the $\sigma$-polynomial, the $w$-polynomial, and the $\tau$-polynomial. This survey examines these relationships, illustrated in Figure~\ref{abst-f1}, and reviews the main results.

\end{abstract}


\begin{figure}[htp] 
	\centering

		\begin{tikzpicture} 
		[auto,node distance=2.3 cm, scale =0.9, transform shape]
		
		\tikzstyle{vertex} = [circle, draw=blue!90]
		
		\tikzstyle{selected vertex} =[rectangle, draw, rounded corners=1.5mm] (rounded) {Rounded};

		\node[selected vertex] (v1)
		at (0,5) {Tutte polynomial};
		
		\node[selected vertex] (v2)
		at (0,3)
		{Characteristic polynomial};

		\draw[-{Latex[length=6mm, width=4mm]}, line width=1.5mm] (v1) -- (v2);
		
		\node[selected vertex] (v31)
		at (-4,1)
		{Chromatic polynomial};
		
		\node[selected vertex] (v32)
		at (4,1) {Flow polynomial};
		
		\draw[-{Latex[length=6mm, width=4mm]}, line width=1.5mm] (-1,2.6) -- (v31);
		
		\draw[-{Latex[length=6mm, width=4mm]}, line width=1.5mm] (1,2.6) -- (v32);
		
		\draw[ -{Latex[length=6mm, width=4mm]}, line width=1.5mm, ->, -Stealth] (-1.5,1) -- (2,1);
		
		\draw[ -{Latex[length=6mm, width=4mm]}, line width=1.5mm, ->, -Stealth]  (1,1) -- (-1.7,1);

		\node[selected vertex] (v41)
		at (-8,-1) {$\sigma$-polynomial};
		
		\node[selected vertex] (v42)
		at (-4.5,-1) {$w$-polynomial};
		
		\node[selected vertex] (v43)
		at (-1,-1) {$\tau$-polynomial};
		
		\draw[-{Latex[length=6mm, width=4mm]}, line width=1.5mm] (-5.5,0.6) -- (v41);

		\draw[-{Latex[length=6mm, width=4mm]}, line width=1.5mm] (-4.5,0.6) -- (v42);
		
		\draw[-{Latex[length=6mm, width=4mm]}, line width=1.5mm] (-3.0,0.6) -- (v43);

	\end{tikzpicture}

	\caption{Relation 
		between chromatic polynomial and other
	polynomials
	\label{abst-f1}
}
\end{figure}

\tableofcontents


\pagestyle{myheadings}

\markboth{Section \thesection}{A survey of  polynomials related to chromatic polynomials}


\section{Tutte polynomial $T_G(x,y)$}
\relabel{sect-tutte}


\subsection{Definition} 


The roots of the Tutte polynomial go back to the work of William Tutte in the 1940s - 1950s. Tutte was interested in graph invariants that could be defined recursively via deletion and contraction of edges
(see \cite{tuttew}). Earlier, similar ideas had appeared in the work of Hassler Whitney, particularly in his study of chromatic polynomials and matroids.

 For a multigraph $G=(V,E)$,
its {\it Tutte polynomial}, denoted by $T_G(x,y)$, 
is defined as follows (see 
Definition 2.3 in \cite{Jo2022}):
\equ{sec2-e1}
{
T_G(x,y)=\sum_{A\subseteq E} (x-1)^{r(E)-r(A)} (y-1)^{|A|-r(A)},
}
where $r(A)=|V|-c(A)$ and $c(A)$ is the number of components of the spanning 
subgraph $(V,A)$.

For any edge $e$ of $G$, 
$G\backslash e$ is defined to be the graph obtained from $G$ by removing $e$, and 
if $e$ is not a loop, 
$G/e$ is defined to be the graph obtained from $G$ 
by contracting $e$
(i.e., the graph obtained from 
$G\backslash e$ by identifying the two ends of $e$).
$T_G(x,y)$
can be obtained recursively 
by deletion
and contraction.

\begin{theo}[Tutte \cite{tut1947}]
	\label{Tutte-recur}
For any multigraph $G=(V,E)$, 
$T_G(x,y)$ can be determined by the following rules:
\begin{enumerate}[(a)]
\item $T_G(x,y)=1$ if $E=\emptyset$; 
\item $T_G(x,y)=yT_{G\backslash e}(x,y)$, where $e$ is a loop of $G$;
\item $T_G(x,y)=xT_{G/e}(x,y)$,
where $e$ is a bridge of $G$; and 
\item $T_G(x,y)=T_{G/e}(x,y)+T_{G\backslash e}(x,y)$, 
where $e$ is neither a bridge 
nor a loop of $G$.
\end{enumerate}
\end{theo}

Let $w$ be an injective weight function $w: E\mapsto \Z$,
where $\Z$ is the set of positive integers. 
For any spanning tree $T$ of $G$ and $e\in E(T)$, $T\backslash e$ has exactly two components,
say $T_1$ and $T_2$.
If $w(e)\le w(e')$ holds for every  $e'\in E(G)$ which has one end 
in $V(T_1)$ and the other end in 
$V(T_2)$, 
$e$ is called an \red{\it internally active edge}
with respect to $T$.
For an edge $e\in E(G)\setminus E(T)$, $e$ is called an \red{\it externally active edge}
with respect to $T$ if $w(e)\le w(e')$ holds for every edge $e'$ on the 
unique cycle of the spanning subgraph $(V, E(T)\cup \{e\})$.

The Tutte polynomial $T_G(x,y)$ 
can also be expressed 
in terms of spanning trees of $G$ (\cite{tut1947,tut1}). 

\begin{theo}[Tutte \cite{tut1947, tut1}] 
	\label{Tutte-tree}
	For any connected multigraph $G=(V,E)$, 
$$
T_G(x,y)=\sum_{T} x^{in(T)} y^{ex(T)},
$$
where the sum runs over all spanning trees $T$ of $G$, 
and  $in(T)$ (resp. $ex(T)$) is the number of internally active edges 
 (resp. externally active edges) with respect to $T$.
 \end{theo}

\begin{exa}\label{Tutte-ex1}
	\begin{enumerate}[(a)]
\item If each edge in $G$ is either a bridge or a loop, then 
$T_G(x,y)=x^{b}y^{l}$,
where $b$ and $l$ are 
the numbers of bridges and 
loops in $G$, respectively.

\item If $G$ is a tree of order $n$, then $T_G(x,y)=x^{n-1}$.

\item If $G$ is a cycle of order $n$, then 
$
T_{G}(x,y)=x+x^2+\cdots+x^{n-1}+y.
$
\end{enumerate} 
\end{exa} 

The conclusions in Example~\ref{Tutte-ex1}
can be obtained by applying  (\ref{sec2-e1}), Theorem~\ref{Tutte-recur}
or Theorem~\ref{Tutte-tree}.

The Tutte Polynomial of 
a graph can be naturally extended to that of a 
matroid.
A matroid $M$ can be 
expressed as $(E,r)$, 
where $E$ is its 
ground set and $r$ is its  
rank function 
from $2^E$ to $\Z_0$
(i.e., the set of non-negative integers), satisfying the following three conditions
(see \cite{oxley2006, oxley2022}): 
\begin{enumerate}[(a)]
	\item $0\le r(A)\le |A|$ for all $A\subseteq E$; 
	\item $r(A)\le r(B)$ for any 
	$A,B\subseteq E$ with  $A\subseteq B$; and 
	\item (submodularity) 
	$r(A\cup B)+r(A\cap B)\le r(A)+r(B)$ 
	for any $A,B\subseteq E$.
\end{enumerate}
The Tutte polynomial of $M=(E,r)$,
denoted by $T_M(x,y)$,  
is defined as follows
(\cite{Crapo1969, oxley2022}):
$$
T_M(x,y)=\sum_{A\subseteq E}(x-1)^{r(E)-r(A)}
(y-1)^{|A|-r(A)}.
$$
For any graph $G$, if $M_G$ 
is the graphic matroid of $G$, 
then $T_G(x,y)=T_{M_G}(x,y)$.

 Theorems~\ref{Tutte-recur} and \ref{Tutte-tree} remain valid when a graph $G$ is replaced by a matroid $M$. In particular, under this substitution in Theorem~\ref{Tutte-tree}, the Tutte polynomial $T_M(x,y)$ is given by the sum of $x^{in(B)}y^{ex(B)}$ over all bases $B$ of $M$ (see \cite{oxley2022}),
 where $in(B)$ is the number 
of internally active members 
$e\in B$
and a member $e\in B$ 
is said to be internally active
with respect to $B$ 
if $e$ is the smallest element (under an injective mapping $w:E\mapsto \Z$) 
in its unique basic bond 
$bo_M(B,e)$,
while $ex(B)$ is the number 
of externally active members 
$e\in E\setminus B$
and an element $e\in E\setminus B$ is said to be 
externally active 
with respect to $B$ 
if $e$ is the smallest element
in the unique fundamental circuit $C(B,e)$.



 \label{spe-case}

\subsection{Connection to other polynomials}

Let $M=(E,r)$ be a matroid, and let $G=(V,E)$ be a graph with $c(G)$ components. 
Relations between Tutte polynomial and some other polynomials are provided below.

\begin{enumerate}[(i)]
	\item {\bf Characteristic polynomial}: 
	
 $C(M,x)=(-1)^{r(M)}
	T_M(1-x,0)$,
	where $C(M, x)$ is the characteristic polynomial of $M$ (see the definition at  (\ref{char-fun})).

	\item  {\bf Chromatic polynomial}: 
	
 $\chi(G,x)=(-1)^{|V|-c(G)}x^{c(G)} T_G(1-x,0)$,
	where $\chi(G,x)$ is the chromatic polynomial of $G$
	(see its definition in Page~\pageref{sect-chro}).
	
	\item  {\bf Flow polynomial}:

	$F(G,y)=
	(-1)^{|E|-|V|+c(G)}T_G(0,1-y) $,
	where $F(G,x)$ is the flow polynomial of $G$ (see its definition at (\ref{eq1-2})).

	\item  {\bf Reliability polynomial}: 
	
	 $R_G(p)=(1-p)^{|E|-|V|+c(G)}
	 p^{|V|-c(G)}
	 T_G(1,\frac 1{1-p})$,
	 where $R_G(p)$ is the reliability polynomial of 
	 a connected $G$,
	 which is defined as follows
	 (see \cite{mon2008, god2001, pag1994, roy2004}):
	 \equ{Tutte-e2}
	 {
	 R_G(p)=\sum_{H\in {\cal G}(G)}
	 p^{|E(H)|}(1-p)^{|E|-|E(H)|},
	}
	 where ${\cal G}(G)$
	 is the set of connected spanning subgraphs of $G$.

\item {\bf Potts Model Partition Function (Statistical Physics)}:

The $q$-state Potts model
(or the multivariable Tutte polynomial) 
of $G$
is defined as follows:
\begin{equation}
	\relabel{pott-eq2}
	Z_G(q,\{w_e\})
	=\sum_{A\subseteq E}q^{c(A)}
	\prod_{e\in A}w_e,
\end{equation}
where $w_e$ is the weight of edge $e$ and 
$c(A)$ is the number of components 
of the spanning subgraph 
$(V,A)$.
Then,
$$
Z_G(q,v)=q^{c(G)}v^{|V|-c(G)}
T_G\Bra{1+\frac qv, v+1},
$$
where $Z_G(q,v)=Z_G(q,\{w_e\})$
with $w_e=v$ for all $e\in E$
(see \cite{sok2005}).
\end{enumerate}

\subsection{Identities} 

\begin{enumerate}[(1)]
\item {\bf Duality:}

\begin{pro}[\cite{tut1}]
\relabel{tutte-pro1}
If $G$ is a connected plane graph and $G^*$ is its dual, then 
$$
T_G(x,y)=T_{G^*}(y,x).
$$
\end{pro}


For any matroid $M=(E,r)$, 
the {\it dual matroid} of $M$, 
denoted by $M^*$, 
is defined to be the one $(E,r^*)$ with its rank function $r^*(A)$ 
determined  by 
$$
r^*(A)=|A|-\big (r(E)-r(E-A)\big ).
$$
for any $A\subseteq E$, i.e., 
$r^*(A)=|A|-\min\limits_{B\in \setb(M)} |B\cap A|$,
where $\setb(M)$ is the family of bases of $M$.

\begin{pro}[\cite{oxley2022}]
\relabel{tutte-pro2}
If $M=(E,r)$ is a matroid and $M^*$ is its dual, then 
$$
T_M(x,y)=T_{M^*}(y,x).
$$
\end{pro}


\item {\bf Factorization}:

\sgap 

If $G$ is disconnected with components $G_1, G_2,\cdots,G_k$, 
or $G$ is connected with blocks $G_1, G_2,\cdots,G_k$,
then 
$$
T_G(x,y)=\prod_{i=1}^k T_{G_i}(x,y).
$$

\item {\bf Coefficients $t_{i,j}$
in $T_M(x,y)$}:

\sgap

\begin{pro}[\cite{lau2011}]\relabel{tutte-pro3}
	Let $M = (E,r)$ be a matroid.
If $t_{i,j}$ is the number of bases $B\in \setb(M)$  
with $in(B)=i$ and $ex(B)=j$,
then 
$$
T_M(x,y)=\sum_{i,j}t_{i,j}x^iy^j,
$$
and the following hold,
where 
$M$ has neither loops nor coloops for (b), (c) and (d):
\begin{enumerate}[(a)]
\item $t_{1,0} =t_{0,1}$  when $|E|\ge 2$;
\item $t_{i,j} =0$ whenever $i > r(M)$ or $j > |E|-r(M)$; 
\item $t_{r(M),0} =1$ and $t_{0,|E|-r(M)} =1$; 
\item $t_{r(M),j}=0$ for all $j > 0$ and
$t_{i,|E|-r(M)}=0$ for all $i > 0$;
and 
\item (\cite{oxley1}) for all $k=0,1,2,\cdots,|E|-1$,
$$
\sum_{i=0}^k\sum_{j=0}^{k-i}(-1)^j{k-i\choose j}t_{i,j}=0
$$
\end{enumerate}
\end{pro}

Note that \red{if $M$ is replaced by a connected graph $G=(V,E)$, then 
Proposition~\ref{tutte-pro3} holds with $r(M)=|V|-1$
and $t_{i,j}$ to be the number of spanning trees $T$ of $G$ 
with $in(T)=i$ and $ex(T)=j$.}

\gap

\item {\bf  An expression
	for $T_M((v+1)/v,v+1)$}:

\begin{theo}
	[\cite{oxley1}]
	\relabel{tut-poly-iden}
For any matroid $M=(E,r)$ and $v\ne 0$, 
$$
T_M((v+1)/v,v+1)=\frac{(v+1)^{|E|}}{v^{r(M)}}.
$$
In particular, for any connected graph $G$ of order $n$ and size 
$m$, 
$$
T_G((v+1)/v,v+1)=\frac{(v+1)^{m}}{v^{n-1}}.
$$
\end{theo}

\proof Let $f_M(u,v)=\sum\limits_{A\subseteq E}u^{r(M)-r(A)}v^{|A|}$.
Then 
$$
f_M(1,v)=(v+1)^{|E|}.
$$
It can be shown that 
$$
f_M(u,v)=v^{r(M)}T_M(u/v+1,v+1),
$$
implying that 
$$
f_M(1,v)=v^{r(M)}T_M(1/v+1,v+1)
=v^{r(M)}T_M((v+1)/v,v+1).
$$
Thus the result follows.
\proofend

By Theorem~\ref{tut-poly-iden}, 
it can be shown that 
$$
(v+1)^{|E|}=\sum_{i,j\ge 0}t_{i,j}v^{r(M)-i}(v+1)^{i+j}.
$$
Let $w=v+1$. Then 
$$
w^{|E|}
=\sum_{i,j\ge 0}t_{i,j}w^{i+j}(w-1)^{r(M)-i}
=\sum_{i,j\ge 0}t_{i,j}\sum_{k=0}^{r(M)-i} 
(-1)^k {r(M)-i\choose k}w^{r(M)+j-k}.
$$
Thus 
$$
w^{|E|-r(M)}=\sum_{i,j\ge 0}t_{i,j}\sum_{k=0}^{r(M)-i} 
(-1)^k {r(M)-i\choose k}w^{j-k}.
$$
If $|E|>r(M)$, then the constant term of the right-hand-side
is $0$, implying that 
$$
\sum_{i=0}^{r(M)}\sum_{j=0}^{r(M)-i} (-1)^j 
{r(M)-i\choose j}t_{i,j}=0.
$$
More identities can be obtained by taking $k$ such that 
$j-k\ne |E|-r(M)$.


\gap

\item 
{\bf A Convolution Formula for $T_M(x,y)$}:

\begin{theo}[Kook et al \cite{kok1999}]
\relabel{tut-theo1}
For any matroid $M=(E,r)$,
\equ{sec-eq5}
{
T_M(x,y)=\sum_{A\subseteq E} T_{M/A}(x,0) T_{M|A}(0,y).
}
\end{theo}

Note that $T_{M/A}(x,0)=0$ if $M/A$ has a loop 
and $T_{M|A}(0,y)=0$ if $M|A$ has a bridge.
Thus (\ref{sec-eq5})
can be revised as follows:

\equ{sec-eq6}
{
T_M(x,y)=\sum_{A\subseteq \F^*(M)} T_{M/A}(x,0) T_{M|A}(0,y).
}
where $\F^*(M)$ is the family of 
those flats $F$ which contain no bridge.

For a graph $G=(V,E)$, let $\P^*(G)$ be the 
family of those partitions $P=\{V_1, V_2,\cdots,V_r\}$ of $V$ 
such that $V_i\ne \emptyset$ and the induced subgraph $G[V_i]$ is connected and bridgeless for each $i$.
For any $P=\{V_1, V_2,\cdots,V_r\}\in \P^*(G)$,
let $G|P$ denote the spanning subgraph of $G$ 
which is the disjoint union of $G[V_i]$'s for $i=1,2,\cdots,r$,
and let $G/P$ be the graph obtained from $G$ by 
contracting all edges in $G[V_i]$ for all $i=1,2,\cdots,r$.
Thus $G/P$ is a graph of order $r$.

 
 Thus, the following conclusion
 follows from (\ref{sec-eq6}).

\begin{theo}\relabel{tut-theo3}
For any graph $G=(V,E)$,
$$
T_G(x,y)=\sum_{P\subseteq \P^*(G)} T_{G/P}(x,0) T_{G|P}(0,y).
$$
\end{theo}

\end{enumerate}

\subsection{Combinatorial Interpretations}
\label{sec1-4}

	

Assume that $G=(V,E)$ is  connected.

\begin{enumerate}[(1)]

\item {\bf $T_G(i,j)$ for $i,j\in \{0,1,2\}$}:
	
\begin{enumerate}[(a)]
\item $T_G(0,0)=0$ if $E\ne \emptyset$;

\item $T_G(2,2)=2^{|E|}$
(i.e., the number of subsets of $E$);

\item $T_G(1,2)$ is the number of spanning 
connected subgraphs of $G$;

\item 
$T_G(2,1)$ is the number of spanning forests of $G$;

\item $T_G(1,1)$  is  the number of spanning trees of $G$,
denoted by $\tau(G)$;

\item $T_G(0,1)$ is the number of those spanning trees $T$ of $G$
with $in(T)=0$;

\item 
$T_G(1,0)$ is the number of those spanning trees $T$ of $G$
with $ex(T)=0$;


\item $T_G(2,0)$ is the number of acyclic orientations of $G$,
denoted by $\alpha(G)$, and 

\item
$T_G(0,2)=$ is 
the number of totally cyclic orientations of $G$,
denoted by $\alpha^*(G)$.
\end{enumerate} 

Properties (a)-(g) above 
follow directly from the definition of  $T_G(x,y)$, 
while (h) is a special case of Theorem~\ref{stanley73}, 
and for 
(i), we refer the readers to
\cite{ver1980}.

\item {\bf $T_G(k+1,0)$ for 
any positive integer $k$}:

\begin{theo}
[Stanley \cite{sta1973a}]
\label{stanley73} 
For any integer $k\ge 1$, $T_G(k+1,0)$ is equal to
$$
\frac 1k |\chi(G,-k)|=\frac 1k\sum_{j=1}^k (k)_j|\Upsilon_j|,
$$
where $\Upsilon _j$ is the set of order pairs
$(P,\O)$,
where $P$ is a partition of $V$ into exactly $j$ non-empty subsets 
$V_1,V_2,\cdots, V_{j}$ and 
$\O$ is an acyclic orientation of the spanning subgraph of $G$ 
with edge set 
$\cup_{1\le i\le j} E(G[V_i])$.
\end{theo}

\proof By Stanley's result in \cite{sta1973a},
$$
|\chi(G,-k)|=\tilde \chi(G,k),
$$
where $\tilde \chi(G,k)$ is the number of order pairs 
$(f,D)$, where $D$ is an acyclic orientation of $G$ 
and $f$ is a mapping $f:V\mapsto \{1,2,\cdots,k\}$ 
such that $f(u)\le f(v)$ whenever $u\rightarrow v$ is an arc in $D$.

For any $j$ with $1\le j\le k$, 
let $\psi_j$ be the number of order pairs 
$(f,D)$ such that 
\begin{enumerate}[(a)]
	\item
$D$ is an acyclic orientation of $G$;
\item $f$ is a mapping $f:V\mapsto \{1,2,\cdots,k\}$ 
such that $f(u)\le f(v)$ whenever $u\rightarrow v$ is an arc 
in $D$;
and

\item $|f(V)|=j$.
\end{enumerate} 

Thus 
$$
\tilde \chi(G,k)=\sum_{j=1}^k \psi_j
=\sum_{j=1}^k (k)_j |\Upsilon_j|.
$$
\proofend

\item $T_G(-1,-1)$: 

\begin{theo} [Read and Rosenstiehl \cite{read1978}]
\label{Read1978} 
$$
T_G(-1,-1)=(-1)^{|E|}(-2)^{\dim (\setb)},
$$
where $\setb=\setc \cap \setc^{\perp}$ and 
$\setc$ is {\it the cycle space} of $G$ over the finite field $GF(2)$
and $\setc^{\perp}$ is the cocycle space of $G$ over the finite field $GF(2)$.
\end{theo}

Let $E(G)=\{e_i:1\le i\le m\}$ be the edge set of $G$,
and $T$ be a spanning tree of $G$. 

For any edge $e\in E(G)-E(T)$,
there is a unique cycle,
denoted by $C(e)$, on the spanning subgraph of $G$ with edge set $E(T)\cup \{e\}$.
Then, define a vector corresponding to $C(e)$:
$$
X_{e}=(x_1,x_2,\cdots,x_m),
$$
where $x_j=1$ if edge $e_j$ is contained in $C(e)$, and $x_j=0$ otherwise.

{\bf The cycle space $\setc$} is the set of linear combinations 
of all vectors $X_{e}$, where $e\in E(G)\setminus E(T)$,  over $GF(2)$.

For any vertex $v\in V(G)$, we define a vector corresponding to
$E(v)=\{e_i\in E(G): e_i \mbox{ has two different ends and } v 
\mbox{ is one of its ends}\}$:
$$
Y_v=(y_1,y_2,\dots,y_m),
$$
where $y_j=1$ if $e_j\in E(v)$, and $y_j=0$ otherwise. 

{\bf The co-cycle space $\setc^{\perp}$} is the set of linear combinations 
of all vectors $Y_{v}$, where $v\in V(G)$,  over $GF(2)$.

Thus, $\setb$ can be obtained by $\setc$ and $\setc^{\perp}$ obtained above.

\begin{exa}
	Verify 
	Theorem~\ref{Read1978}
	for the graph $G$ below.
	
	\begin{figure}[htp] 
		\centering
		
	
	\begin{tikzpicture}[
		vertex/.style={circle,fill=black,inner sep=0pt,minimum size=7pt},
		every label/.style={font=\small,inner sep=1pt},
		line cap=round,
		line join=round
		]
		\coordinate (v1) at ( 0,     1.15);
		\coordinate (v2) at ( 1.55,  0);
		\coordinate (v3) at ( 0,    -1.25);
		\coordinate (v4) at (-1.55,  0);
		
		\draw[line width=0.9pt] (v4) -- node[above left]  {$e_1$} (v1);
		\draw[line width=0.9pt] (v1) -- node[above right] {$e_2$} (v2);
		\draw[line width=0.9pt] (v2) -- node[below right] {$e_3$} (v3);
		\draw[line width=0.5pt] (v3) -- node[below left]  {$e_4$} (v4);
		\draw[line width=0.5pt] (v4) -- node[above]       {$e_5$} (v2);
		
		\node[vertex,label=above:{$v_1$}]       at (v1) {};
		\node[vertex,label=right:{$v_2$}]       at (v2) {};
		\node[vertex,label=below right:{$v_3$}] at (v3) {};
		\node[vertex,label=above left:{$v_4$}]  at (v4) {};
	\end{tikzpicture}
	
	\caption{A graph $G$
	\label{f2}
} 
	
	\end{figure}

\end{exa}

\solu Note that
$$
T_G(x,y)=x^3 + 2x^2 + 2xy + y^2 + x + y.
$$ 
Thus, $T_G(-1,-1)=2$. 
It remains to show that $\dim (\setb)=1$.

Let $T$ be the spanning tree of $G$ with edge set $\{e_1,e_2,e_3\}$. 
Then, $\setc$ is the linear space over $GF(2)$ spanned by vectors: 
$$
(1,1,1,1,0), \quad (1,1,1,0,0).
$$
While $\setc^{\perp}$ is the linear space over $GF(2)$ spanned by vectors: 
$$
(1,1,0,0,0),  (0,1,1,0,1), (0,0,1,1,0),
(1,0,0,1,1).
$$
It is clear that $(1,1,1,1,0)\in \setc^{\perp}$, but $(1,1,1,0,0)\not\in \setc^{\perp}$.

Thus, $\dim(\setb)=1$, and the result is verified. 
\proofend 

\item  {\bf $T_G(1-n,1-n)$, 
$T_G(1+n,1+n)$, 
$T_G(-1,-1)$
and $T_G(3,3)$
for a plane graph $G$,
where $n$ a positive integer}:



Let $G$ be a plane graph. 
The {\it medial graph} of $G$ is constructed by placing a vertex on each edge of $G$ and drawing edges around the faces of $G$. 

The faces of this medial graph are colored black or white, depending on whether they contain or do not contain, respectively, a vertex of the original graph $G$. 

{The directed medial graph of $G$}, denoted by $\tilde G_m$,
is obtained by assigning a diection to each edge of the medial graph 
so that the black face is on the left. 
An example of $\tilde G_m$ is shown in Figure~\ref{f21}.

\begin{figure}[htbp]
  \centering

\begin{tikzpicture}[
	line width=0.35pt,
	line cap=round,
	line join=round,
	original vertex/.style={circle,fill=black,draw=black,
		inner sep=0pt,minimum size=7pt},
	medial vertex/.style={circle,fill=black,draw=black,
		inner sep=0pt,minimum size=3.1pt},
	open vertex/.style={circle,fill=white,draw=black!65,
		line width=0.3pt,inner sep=0pt,minimum size=4.5pt},
	mid arrow/.style={
		decoration={markings,
			mark=at position #1 with {\arrow{Stealth[length=2.2pt,width=1.7pt]}}},
		postaction={decorate}
	}
	]
	
	\begin{scope}
		\coordinate (T) at ( 0,     1.08);
		\coordinate (R) at ( 1.45,  0);
		\coordinate (B) at ( 0,    -1.08);
		\coordinate (L) at (-1.45,  0);
		
		\draw (T) -- (R) -- (B) -- (L) -- cycle;
		\draw (L) -- (R);
		\foreach \v in {T,R,B,L}
		\node[original vertex] at (\v) {};
		
		\node at (0,-1.98) {$G$};
	\end{scope}
	
	\begin{scope}[xshift=4.25cm]
		\coordinate (T)  at ( 0,      1.08);
		\coordinate (R)  at ( 1.45,   0);
		\coordinate (B)  at ( 0,     -1.08);
		\coordinate (L)  at (-1.45,   0);
		\coordinate (NW) at (-0.725,  0.54);
		\coordinate (NE) at ( 0.725,  0.54);
		\coordinate (SW) at (-0.725, -0.54);
		\coordinate (SE) at ( 0.725, -0.54);
		\coordinate (O)  at ( 0,      0);
		
		\draw[black!65,line width=0.25pt]
		(T) -- (R) -- (B) -- (L) -- cycle;
		\draw[black!65,line width=0.25pt] (L) -- (R);
		\foreach \v in {T,R,B,L}
		\node[open vertex] at (\v) {};
		
		\draw[mid arrow=0.77] (NE)
		.. controls (0.20,1.95) and (-0.20,1.95) .. (NW);
		\draw[mid arrow=0.54] (NW)
		.. controls (-0.17,0.90) and (0.17,0.90) .. (NE);
		
		\draw[mid arrow=0.77] (SW)
		.. controls (-0.20,-1.95) and (0.20,-1.95) .. (SE);
		\draw[mid arrow=0.38] (SE)
		.. controls (0.17,-0.90) and (-0.17,-0.90) .. (SW);
		
		\draw[mid arrow=0.91] (NW)
		.. controls (-2.34,0.27) and (-2.34,-0.27) .. (SW);
		\draw[mid arrow=0.94] (SE)
		.. controls (2.34,-0.27) and (2.34,0.27) .. (NE);
		
		\draw[mid arrow=0.51] (O)
		.. controls (-1.10,0) and (-0.98,0.12) .. (NW);
		\draw[mid arrow=0.54] (NE)
		.. controls (0.98,0.12) and (1.10,0) .. (O);
		\draw[mid arrow=0.20] (SW)
		.. controls (-0.98,-0.12) and (-1.10,0) .. (O);
		\draw[mid arrow=0.51] (O)
		.. controls (1.10,0) and (0.98,-0.12) .. (SE);
		
		\foreach \v in {NW,NE,SW,SE,O}
		\node[medial vertex] at (\v) {};
	\end{scope}
	
	\begin{scope}[xshift=8.60cm]
		\coordinate (NW) at (-0.725,  0.54);
		\coordinate (NE) at ( 0.725,  0.54);
		\coordinate (SW) at (-0.725, -0.54);
		\coordinate (SE) at ( 0.725, -0.54);
		\coordinate (O)  at ( 0,      0);
		
		\draw[mid arrow=0.85] (NE)
		.. controls (0.20,1.14) and (-0.20,1.14) .. (NW);
		\draw[mid arrow=0.54] (NW)
		.. controls (-0.17,0.90) and (0.17,0.90) .. (NE);
		
		\draw[mid arrow=0.87] (SW)
		.. controls (-0.20,-1.30) and (0.20,-1.30) .. (SE);
		\draw[mid arrow=0.38] (SE)
		.. controls (0.17,-0.90) and (-0.17,-0.90) .. (SW);
		
		\draw[mid arrow=0.91] (NW)
		.. controls (-1.73,0.15) and (-1.73,-0.22) .. (SW);
		\draw[mid arrow=0.94] (SE)
		.. controls (1.83,-0.25) and (1.83,0.13) .. (NE);
		
		\draw[mid arrow=0.51] (O)
		.. controls (-1.10,0) and (-0.98,0.12) .. (NW);
		\draw[mid arrow=0.54] (NE)
		.. controls (0.98,0.12) and (1.10,0) .. (O);
		\draw[mid arrow=0.20] (SW)
		.. controls (-0.98,-0.12) and (-1.10,0) .. (O);
		\draw[mid arrow=0.51] (O)
		.. controls (1.10,0) and (0.98,-0.12) .. (SE);
		
		\foreach \v in {NW,NE,SW,SE,O}
		\node[medial vertex] at (\v) {};
		
		\node at (0,-1.98) {$\widetilde{G}_m$};
	\end{scope}
	
\end{tikzpicture}

  \caption{The directed medial graph $\tilde G_m$ of $G$
\relabel{f21}
}
\end{figure}

For any directed graph $H$, let 
$\D_n(H)$ be  the family of ordered partitions $(D_1,\cdots,D_n)$
of $E(H)$ such that  $H$ restricted to $D_i$ is
2-regular and consistently oriented for all $i$.

\begin{theo}[Martin \cite{martin1977}] \relabel{tutte-theo10}
Let $\tilde G_m$ be  the directed medial graph of a plane graph $G$. 
 Then, for any positive integer $n$, 
$$
(-n)^{c(G)}T_G(1-n,1-n) = \sum_{(D_1,\cdots,D_n) \in \D_n(\tilde G_m)}
(-1)^{\sum_{1\le i\le n}c(D_i)}
$$ 
and
$$
n^{c(G)}T_G(1+ n,1+ n) =\sum_{\phi}2^{\mu(\phi)}, 
$$
where the sum runs over all edge colorings $\phi$ of $\tilde G_m$  
with $n$ colors so that each (possibly empty) set of monochromatic edges 
forms an Eulerian digraph, and where $\mu(\phi)$ 
is the number of monochromatic vertices in the coloring $\phi$.
\end{theo}

An {\it anticircuit} in a digraph is a closed trail so that the directions of 
the edges alternate as the trail passes through any vertex of degree greater 
than 2.
In 4-regular Eulerian digraph, a anticircuit can be obtained
by choosing the two incoming edges or the two outgoing edges 
at each vertex. 

\begin{theo}\relabel{tutte-theo12}
Let $G$ be a connected plane graph, and 
let $a(\tilde G_m)$ denote
 the number of anticircuits in $\tilde G_m$. 
Then 
\begin{enumerate}[(a)]
\item (Martin \cite{martin1978})
$
T_G(-1,-1) = (-1)^{|E(G)|}(-2)^{a(\tilde G_m)-1} 
$, and 

\item (Vergnas \cite{ver1988})
$
T_G(3,3) = K 2^{a(\tilde G_m)-1}, 
$
where $K$ is some odd integer.
\end{enumerate} 
\end{theo}
 
Note that $a(\tilde G_m)$ is actually equal to 
the number of components of the link diagram $D(G)$.

\item {\bf Expressions 
	for $[y^j]T_G(1,y)$, 
	the coefficients of $y^j$ in $T_G(1,y)$}:

Let $G$ be a $(k + 1)$-edge connected graph of order 
$n$ and size $m$, and 
let $g = m -n + 1$. 
In 2025,  
Guan, Jin
and K\'alm\'an
 \cite{Guan2025} found the following expression for 
$[y^j]T_G(1,y)$:
\equ{guan2025}
{
\forall j: g-k\le j\le g, \quad 
[y^j]T_G(1,y)={m-j-1\choose n-2}.
}
The above result was 
further extended 
to the following one
by Chen and Guo \cite{Chen2025}:
\equ{chen2025}
{
	[y^j]T_G(1,y)
	={m-j-1\choose n-2}
	-\sum_{i=k+1}^{g-j}
	{m-j-i-1\choose n-2}|\EC_i(G)|,
}
for all $ j$ with
$ g-3(k+1)/2< j\le g$,
where $\EC_i(G)$ is the set of minimal edge-cuts 
of size $i$ in $G$.
The results in (\ref{guan2025})
and (\ref{chen2025})
were recently generalized to 
$T_M(x,y)$ for a  matroid
$M=(E, r)$ by Ma, Guan and Jin \cite{Ma2025}:
\equ{ma2025}
{
	[y^j]T_M(1,y)
	=\sum_{t=j}^{|E|-r(M)}
	(-1)^{t-j}{t\choose j}\sigma_{r(M)+t}(M),
}
where  $\sigma_s(M)$ is the number of spanning sets $S$
of $M$ with $|S|=s$.

\end{enumerate} 


\subsection{Invariants of Tutte polynomial} 

	

\begin{pro}
	[\cite{Bry1980}]
	\relabel{tut-pro5}
	For connected matroids $M_1$ and $M_2$,
	if $T_{M_1}(x,y)=T_{M_2}(x,y)$, then 
	\begin{enumerate}[(a)]
		\item $r(M_1)=r(M_2)$;
		\item $|E(M_1)|=|E(M_2)|$;
		\item for each $i$ with $0\le i \le r(M_1)$, 
		the number of independent sets of $M_1$ of cardinality $i$
		is equal to the number of independent sets of $M_2$ of cardinality $i$;
		\item the girth $g(M_1)=g(M_2)$;
		\item the number of circuits of $M_1$ of cardinality $g(M_1)$
		is equal to the number of circuits of $M_2$ of cardinality $g(M_2)$;
		\item for each $i$ with $0 \le i\le r(M_1)$, 
		if $f_i(M_j)$ is the largest cardinality among all flats of $M_j$ of rank $i$, 
		then $f_i(M_1)=f_i(M_2)$;
		and 
		
		\item the number of rank-$i$ flats $F_1$ of $M_1$ with $|F_1|=f_i(M_1)$
		is equal to the number of rank-$i$ flats $F_2$ of $M_2$ with $|F_2|=f_i(M_2)$.
	\end{enumerate}
\end{pro}

\begin{pro}\relabel{tut-pro6}
	Let $G$ be a simple $2$-connected graph. 
	The following parameters of $G$  are determined by its Tutte polynomial 
	$T_G(x,y)$:
	\begin{enumerate}[(a)]
		\item  the edge-connectivity $\lambda(G)$; 
		in particular, a lower bound for the minimum degree $\delta(G)$; 
		\item the number of cliques of each order and the clique-number $\omega(G)$;  and 
		\item the number of cycles of length three, four and five, and 
		the number of cycles of length four with exactly one chord.
	\end{enumerate}
\end{pro}

\subsection{Universality of the Tutte Polynomial}

\begin{theo}[Brylawski and Oxley~\relabel{oxley1}]
\relabel{tutte-theo4}
Let $\G$ be a minor closed class of graphs. 
If a graph invariant $f$ from $\G$ to a commutative ring $R$ with unity 
satisfying all conditions below for any $G,H\in \G$:
\begin{enumerate}[(a)]
\item $f(N_1)=1$, where $N_k$ 
is the empty graph of order $k$;
\item $f(G\cup H)=f(G)f(H)$,
where $G\cup H$ is the union 
of vertex-disjoint graphs $G$ and $H$; 
\item $f(G)=x_0f(G/e)$ if $e$ is a bridge;
\item $f(G)=y_0f(G\backslash e)$ if $e$ is a loop; and 
\item $f(G) = af(G\backslash e)+ bf(G/e)$
for each edge $e$ which is not a loop nor a bridge,
where $a,b$ are non-zero constants; 
\end{enumerate}
then 
$$
f(G)=a^{|E(G)|-r(E(G))}b^{r(E(G))}T_G\left (\frac {x_0}b,\frac {y_0}a\right ).
$$
\end{theo}

\subsection{$T$-equivalent graphs
} 

	

Two graphs $G_1$ and $G_2$ having the same Tutte polynomial
are called {\it codichromatic graphs} by Tutte~\cite{tutte}
and also called {\it $T$-equivalent graphs},
denoted by $G\tsim H$.

\begin{theo}[\cite{good2011}]\relabel{tutt-theo3}
	If $H\tsim G$
	and 
	$G$ is a simple outerplanar graph, 
	then $H$ is also outerplanar.  
\end{theo}

It is trivial that two isomorphic graphs 
are $T$-equivalent. 
If two non-isomorphic graphs have isomorphic cyclic matrods, then 
they are also $T$-equivalent~\cite{oxley1}.

A well-known operation for constructing 
such a pair of graphs 
is {\it the Whitney twist}~\cite{whit} 
which changes a graph to another one
by flipping a subgraph 
 at a vertex-cut of size $2$.
 Such two graphs $G$ and $G'$
 are shown in Figure~\ref{f1-0}, where $\{u_1,u_2\}$ 
 is the cut-set chosen from $G$.
 
 \begin{figure}[htbp]
  \centering

\begin{tikzpicture}[scale=0.9,
	line width=0.35pt,
	line cap=round,
	line join=round,
	vertex/.style={circle,draw=black,fill=black,
		inner sep=0pt,minimum size=9pt},
	vertex label/.style={font=\large,inner sep=1pt},
	graph label/.style={font=\Large}
	]
	
	\begin{scope}
		\coordinate (L)  at (0,0);
		\coordinate (B)  at (1.20,0);
		\coordinate (u1) at (2.60,1.00);
		\coordinate (u2) at (2.60,-1.00);
		\coordinate (T)  at (3.88,1.67);
		\coordinate (R)  at (5.00,0.025);
		\coordinate (D)  at (3.88,-1.30);
		\coordinate (Q)  at ($(L)!0.555!(u1)$);
		
		\draw (u1) -- (T) -- (R) -- (D) -- (u2);
		\draw (u1) -- (R);
		\draw (u2) -- (T);
		
		\draw (L) -- (Q) -- (u1);
		\draw (Q) -- (B);
		\draw (u1) -- (B) -- (u2);
		
		\draw (L) -- (u2);
		\draw (L) .. controls (0.70,-0.72) and (1.75,-0.95) .. (u2);
		
		\foreach \v in {L,B,Q,u1,u2,T,R,D}
		\node[vertex] at (\v) {};
		
		\node[vertex label,anchor=south west]
		at ($(u1)+(-0.04,0.23)$) {$u_1$};
		\node[vertex label,anchor=north west]
		at ($(u2)+(-0.04,-0.14)$) {$u_2$};
		\node[graph label] at (2.80,-1.98) {$G$};
	\end{scope}
	
	\begin{scope}[xshift=6.40cm]
		\coordinate (L)  at (0,0);
		\coordinate (B)  at (1.20,0);
		\coordinate (u1) at (2.60,1.00);
		\coordinate (u2) at (2.60,-1.00);
		\coordinate (T)  at (3.88,1.67);
		\coordinate (R)  at (5.00,0.025);
		\coordinate (D)  at (3.88,-1.30);
		\coordinate (Q)  at ($(L)!0.555!(u2)$);
		
		\draw (u1) -- (T) -- (R) -- (D) -- (u2);
		\draw (u1) -- (R);
		\draw (u2) -- (T);
		
		\draw (L) -- (Q) -- (u2);
		\draw (Q) -- (B);
		\draw (u1) -- (B) -- (u2);
		
		\draw (L) -- (u1);
		\draw (L) .. controls (0.70,0.80) and (1.75,0.99) .. (u1);
		
		\foreach \v in {L,B,Q,u1,u2,T,R,D}
		\node[vertex] at (\v) {};
		
		\node[vertex label,anchor=south west]
		at ($(u1)+(-0.04,0.23)$) {$u_1$};
		\node[vertex label,anchor=north west]
		at ($(u2)+(-0.04,-0.14)$) {$u_2$};
		\node[graph label] at (2.10,-1.98) {$G'$};
	\end{scope}
	
\end{tikzpicture}


  \caption{$G'$ is obtained from $G$ by a Whitney twist.
\relabel{f1-0}
}
\end{figure}


Mentioned by Tutte~\cite{tutte}, the two graphs  
$G_0$ and $H_0$ in  Figure~\ref{f-G10} 
were found by Dr. Marion C. Gray in 1930s.
These two graphs are not isomorphic 
and even have non-isomorphic cyclic matroids,
because $H_0$, unlike $G_0$, contains 
a triangle having no common edge 
with any other triangle \cite{tutte}. 
However, $G_0\backslash e\cong H_0\backslash f$ and 
$G_0/e\cong H_0/f$,
where $e$ and $f$ are the edges in $G_0$ and $H_0$ 
which are expressed by dashed lines in Figure~\ref{f-G10}.
Hence $G_0\tsim H_0$
by Theorem~\ref{Tutte-recur} (d).

\begin{figure}[h!]
  \centering

\begin{tikzpicture}[
	line width=0.35pt,
	line cap=round,
	line join=round,
	vertex/.style={circle,draw=black,fill=black,
		inner sep=0pt,minimum size=8pt},
	marked edge/.style={dash pattern=on 3.5pt off 3pt},
	edge label/.style={font=\huge,inner sep=2pt},
	graph label/.style={font=\huge}
	]
	
	\begin{scope}
		\coordinate (A) at (0,2.40);
		\coordinate (B) at (2.40,2.40);
		\coordinate (D) at (0,0);
		\coordinate (E) at (2.40,0);
		\coordinate (T) at (1.18,3.60);
		\coordinate (C) at (1.25,1.15);
		
		\draw (D) -- (A) -- (B) -- (E);
		\draw[line width=0.85pt] (D) -- (E);
		\draw (A) -- (T) -- (B);
		
		\draw (A) .. controls (0.28,2.99) and (0.75,3.40) .. (T);
		
		\draw (A) -- (C) -- (D);
		\draw[marked edge] (C)
		-- node[edge label,midway,above left] {$e$} (B);
		
		\foreach \v in {A,B,C,D,E,T}
		\node[vertex] at (\v) {};
		
	\end{scope}
	
	\begin{scope}[xshift=6.60cm]
		\coordinate (A) at (0,2.40);
		\coordinate (B) at (2.40,2.40);
		\coordinate (D) at (0,0);
		\coordinate (E) at (2.40,0);
		\coordinate (T) at (1.18,3.60);
		\coordinate (C) at (1.25,1.15);
		
		\draw (A) -- (B) -- (E) -- (D);
		\draw[line width=0.85pt] (D) -- (A);
		\draw (A) -- (T) -- (B);
		
		\draw (A) -- (C) -- (D);
		\draw (A) .. controls (0.22,1.77) and (0.68,1.31) .. (C);
		
		\draw[marked edge] (C)
		-- node[edge label,midway,above right] {$f$} (E);
		
		\foreach \v in {A,B,C,D,E,T}
		\node[vertex] at (\v) {};
		
	\end{scope}
	
\end{tikzpicture}

(a) $G_0$ \hspace{5 cm} (b) $H_0$

 \caption{$G_0\backslash e\cong H_0\backslash f$ and 
$G_0/e\cong H_0/f$
\relabel{f-G10}
}
\end{figure}

Let $\Phi$ be the set of quaternions $(G,e,H,f)$, 
where $G$ and $H$ are loopless graphs
and $e$ and $f$ are edges in $G$ and $H$, respectively, 
such that both 
$G\backslash e\tsim H\backslash f$ 
and $G/e\tsim H/f$ hold.
Clearly, 
$G\tsim H$ holds
for every $(G,e,H,f)\in \Phi$
by Theorem~\ref{Tutte-recur}.
The following problem was proposed and studied in~\cite{dong2025a}.

\begin{prob}
	[\cite{dong2025a}]
	\label{prob2}
	Searching for members $(G,e,H,f)$ of $\Phi$.
\end{prob}

Note that $\Phi$ contains a subfamily $\Phi'$ 
which consists of those quaternions $(G,e,H,f)\in \Phi$ 
such that $G\backslash e\cong H\backslash f$
and $G/e\cong H/f$.
Clearly, $(G_0,e,H_0,f)\in \Phi$
for $G_0$ and $H_0$ shown in Figure~\ref{f-G10}.
In \cite{dong2025a}, 
the authors present a method 
for generating new 
members in $\Phi'$
from any known member $(G,e,H,f)$ in $\Phi'$.

\gap

\begin{enumerate}[(1)]
\item 
{\bf $T$-equivalent graphs produced by flipping a rotor}:

\sgap 
 
Assume that $R$ is a graph and
$\psi$ is an automorphism of $R$.
For any vertex $x$ in $R$, the set  
$\{\psi^i(x): i\ge 0\}$ is called a {\it vertex orbit}
of $\psi$ and $x$ is called a {\it fixed vertex} of $\psi$
if $\psi(x)=x$.

If $R$ is a subgraph of a graph $G$,
a subset $B$ of $V(R)$
is called a {\it border} of $R$ in $G$ if
every edge in $G$ incident with some vertex  
$V(R)-B$ must be an edge in $R$.
We call $R$ a {\it rotor} of $G$ with a border 
$B$ if $B$ is a vertex orbit of 
some automorphism $\psi$.

Tutte~\cite{tutte} showed that if 
$G$ is a graph containing a rotor $R$   
with a border $B$ of size at most $5$,  
then  $G$ and $G'$ are $T$-equivalent, 
where $G'$ is the graph 
obtained from $G$ by flipping $R$ along 
its border $B$, i.e., 
by replacing $R$ by its mirror image. 
We will express Tutte's result below.

Given any vertex-disjoint graphs $G$ and 
$W$ with $\{u_1,u_2,\cdots,u_k\}\subseteq V(G)$ 
and $\{w_1,w_2,\cdots,w_k\} \subseteq V(W)$,
let $G(u_1,u_2,\cdots,u_k)\sqcup 
W(w_1,w_2,\cdots,w_k)$ denote the graph obtained 
from $G$ and $W$ by identifying 
$u_i$ and $w_i$ as a new vertex 
for all $i=1,2,\cdots,k$.
An example of $G(u_1,u_2,u_3)\sqcup
W(w_1,w_2,w_3)$ is shown in Figure \ref{f1}. 

\begin{figure}[htbp]
  \centering
  

\newcommand{\SetVertices}{%
	\coordinate (T) at (0,3.90);
	\coordinate (M) at (0,1.75);
	\coordinate (B) at (0,0);
}

\def\LeftRegion{%
	(T)
	.. controls (-3.05,2.05) and (-3.00,0.90) .. (B)
	.. controls (-1.09,0.63) and (-1.10,1.23) .. (M)
	.. controls (-1.03,2.40) and (-1.02,3.04) .. (T)
	-- cycle
}

\def\RightRegion{%
	(T)
	.. controls (3.90,2.56) and (3.32,1.10) .. (B)
	.. controls (1.10,0.75) and (1.02,1.29) .. (M)
	.. controls (0.84,2.42) and (0.89,3.13) .. (T)
	-- cycle
}

\newcommand{\DrawLeftRegion}{%
	\begin{scope}
		\clip \LeftRegion;
		\foreach \k in {-2,...,13} {
			\draw[hatching]
			(-3.5,{-3.15+0.46*\k}) -- (0.4,{0.36+0.46*\k});
		}
	\end{scope}
	\draw[boundary] \LeftRegion;
}

\newcommand{\DrawRightRegion}{%
	\begin{scope}
		\clip \RightRegion;
		\foreach \k in {-4,...,12} {
			\draw[hatching]
			(-0.4,{-0.22+0.46*\k}) -- (3.6,{1.98+0.46*\k});
		}
	\end{scope}
	\draw[boundary] \RightRegion;
}

\newcommand{\DrawVertices}{%
	\foreach \v in {T,M,B}
	\node[vertex] at (\v) {};
}

	\begin{tikzpicture}[
		line cap=round,
		line join=round,
		boundary/.style={line width=0.35pt},
		hatching/.style={line width=0.25pt,dash pattern=on 2.8pt off 2.5pt},
		vertex/.style={circle,draw=black,fill=black,
			inner sep=0pt,minimum size=8.5pt},
		vertex label/.style={font=\Large,inner sep=1.5pt},
		panel label/.style={font=\Large}
		]
		
		\begin{scope}
			\SetVertices
			\DrawLeftRegion
			\DrawVertices
			
			\node[vertex label,anchor=west] at (0.22,3.95) {$u_3$};
			\node[vertex label,anchor=west] at (0.22,1.75) {$u_2$};
			\node[vertex label,anchor=west] at (0.22,-0.05) {$u_1$};
			\node[panel label] at (-0.70,-0.75) {$(\mathrm{a})\ G$};
		\end{scope}
		
		\begin{scope}[xshift=2.65cm]
			\SetVertices
			\DrawRightRegion
			\DrawVertices
			
			\node[vertex label,anchor=east] at (-0.22,3.85) {$w_3$};
			\node[vertex label,anchor=east] at (-0.22,1.70) {$w_2$};
			\node[vertex label,anchor=east] at (-0.22,0.00) {$w_1$};
			\node[panel label] at (0.80,-0.75) {$(\mathrm{b})\ W$};
		\end{scope}
		
		\begin{scope}[xshift=9.10cm]
			\SetVertices
			\DrawLeftRegion
			\DrawRightRegion
			\DrawVertices
			
			\node[panel label] at (-0.10,-0.75) {$(\mathrm{c})$};
		\end{scope}
		
	\end{tikzpicture}

  \caption{Graph $G(u_1,u_2,u_3)\sqcup
W(w_1,w_2,w_3)$
\relabel{f1}
}
\end{figure}

\begin{theo}[Tutte~\cite{tutte}]\relabel{theo-tutte}
Let $R$ be a connected graph with an automorphism $\psi$.
If $\{u_1, u_2, \cdots,u_k\}$
is a vertex orbit of $\psi$
(i.e., $\psi(u_i)=u_{i+1}$ for all $i=1,2,\cdots,k$), where $k\le 5$, then  
the two graphs 
$R(u_{1},\cdots,u_k)\sqcup W(w_1,\cdots,w_k)$ 
and $R(u_{k},\cdots,u_1)\sqcup W(w_1,\cdots,w_k)$
are $T$-equivalent
for an arbitrary graph $W$,
where $w_1,\cdots,w_k$ are distinct vertices in $W$.
\end{theo}

In \cite{dong2025a}, 
Theorem~\ref{theo-tutte} 
was extended to a rotor $R$ of order $k\ge 6$ 
when $W$ satisfies certain conditions.  

\gap

\item {\bf $T$-unique graphs}: 
\sgap

A graph $G$ is said to be {\it T-unique} if 
for any graph $H$, $H\cong G$ whenever $H\tsim G$.

\begin{theo}[\cite{vinu2003}]
\relabel{tutt-theo8}
The following graphs are T-unique:
\begin{enumerate}[(a)]
\item for every set of positive integers $p_1,p_2,\cdots,p_k$, the complete multipartite graph $K_{p_1,p_2,\cdots,p_k}$ is T-unique, 
with the only exception of $K_{1,p}$;
\item $C_n^2$, where $n\ge 3$ and $C^2_n$ is obtained from the cycle graph $C_n$ 
by adding edges joining any two vertices in $C_n$ 
with distance $2$;
\item graph $C_n \times K_2$;
\item the M\"obius ladder $M_n$, where $n\ge 2$, 
which is constructed from an even cycle $C_{2n}$ 
by joining every pair of vertices at distance $n$; and 
\item The $n$-cube $Q_n$, $n\ge 2$, 
which is defined as the product of $n$ copies of $K_2$.
\end{enumerate}
\end{theo}
\end{enumerate} 

\subsection{Inequalities and equalities}

\begin{enumerate}[(1)]
\item {\bf Inequalities}:

\begin{theo}[Merino et al \cite{meri2009}]\relabel{tutt-theo4}
If a matroid $M$ has neither loops nor isthmuses, then 
$$
\max\{T_M(4,0),T_G(0,4)\}\ge T_M(2,2).
$$
\end{theo}

It can be proved by 
applying the fact that $T_M(2,2)=2^{|E|}$ and $T_M(4,0)\ge 4^{r(M)}$ and
$T_M(0,4)\ge 4^{|E|-r(M)}$.

\begin{theo}[Merino et al \cite{meri2009}]\relabel{tutt-theo5}
If a matroid $M=(E,r)$ contains two disjoint bases, then 
$$
T_M(0,2a)\ge T_M(a,a),
$$
for all $a\ge 2$.
Dually, if its ground set $E$ is the union of two bases of $M$, 
then 
$$
T_M(2a,0)\ge T_M(a,a),
$$
for all $a\ge 2$.
\end{theo}


\begin{con}[Merino and Welsh~\cite{meri1999}]
\relabel{tutt-con1}
Let $G$ be a $2$-connected graph without loops. Then 
$$
\max\{T_G(2,0),T_G(0,2)\}\ge T_G(1,1).
$$
\end{con}

Merino and Welsh also mentioned the following conjecture
which strengthens Conjecture~\ref{tutt-con1}.

\begin{con}[Merino and Welsh~\cite{meri1999}]
\relabel{tutt-con2}
Let $G$ be a 2-connected graph without loops. Then 
$$
T_G(2,0)T_G(0,2)\ge T_G(1,1)^2.
$$
\end{con}

Some partial solutions of Conjecture~\ref{tutt-con1}
have been obtained. 
These conjectures have been naturally extended to any matroid $M$. 
However, 
the matroidal
version of this conjecture has been disproved by 
Beke, Cs\'aji,
Csikv\'ari, and Pituk \cite{beke2024} in 2024
who constructed
explicit counterexamples.
Thus,  the conjecture remains open only for graphic matroids (i.e., graphs).

\begin{theo}[Thomassen \cite{tho2010}]
	\relabel{tut-theo6}
Let $G$ be a graph of order $n$ and size $m$.
\begin{enumerate}[(a)]
	\item If G is simple and $m\le 16n/15$ edges, then 
$T_G(2,0) > T_G(1,1)$; and 
\item if $G$ is a bridgeless  and 
$m\ge 4n-4$ edges, then 
$T_G(0,2)> T_G(1,1)$.
\end{enumerate} 
\end{theo}

A graph is called a {\it series-parallel graph}
if it is obtained 
from a single edge by repeatedly duplicating or subdividing edges in any fashion.

\begin{theo}[Noble and Royle \cite{nob2014}]\relabel{tut-theo5}
Conjecture~\ref{tutt-con1} holds for all series-parallel graphs.
\end{theo}

\begin{theo}[Jackson \cite{jac2010}]\relabel{tut-theo7}
Let $G$ be a graph without loops or bridges and 
$a,b$ be positive real numbers with $b \ge a(a + 2)$. 
Then 
$
\max\{T_G(b,0),T_G(0,b)\}\ge T_G(a,a).
$
\end{theo}

\item {\bf Equalities}: 

\begin{theo}[Merino \cite{meri2008}]
$
T_{K_{n+2}}(1, -1)=T_{K_n}(2 , -1).
$
\end{theo} 

Merino's result was generalized by 
Goodall et al \cite{good2009}.
A graph is called a {\it threshold graph} if the vertices can be ordered so that each vertex is adjacent to either all or none of the previous vertices. 
Threshold graphs are also the graphs with no induced $P_4,C_4$ or $2P_2$. 

\begin{theo}[Goodall et al \cite{good2009}]
If $G$ is a threshold graph and $u$ and $v$ are the first and last vertex in an ordering of the vertices of $G$ such that each vertex is adjacent to either
all or none of the previous ones, then 
$T_G(1 , -1)=T_{G-u-v}(2 , -1)$.
\end{theo} 
\end{enumerate}

\subsection{Open problems}

Except for Problem~\ref{prob2} and Conjectures~\ref{tutt-con1} and~\ref{tutt-con2}, we conclude this section by presenting the following problems on Tutte polynomials.

\begin{prob}\label{prob7}
	Except known results in 
	Subsection~\ref{sec1-4}, 
which evaluations $T_G(a,b)$
admit natural combinatorial interpretations?
\end{prob}

\begin{prob}\label{prob4} 
	Characterize all polynomials that arise as Tutte polynomials of graphs or matroids.
\end{prob}

\begin{prob}\label{prob5} 
Are coefficients of 
Tutte polynomials unimodal for broad classes of graphs or matroids?
\end{prob}

\begin{prob}\label{prob8} 
Given any graph $G$, determine all graphs $H$ 
such that $H\tsim G$.
\end{prob}


\section{Characteristic polynomial $C(M,x)$}
\relabel{sect-char}

\subsection{Definition}


The concept of characteristic polynomial of a matroid 
was  introduced by Rota 
\cite{rot1964}.
The {\it characteristic polynomial}
$C(M, x)$ of a matroid $M=(E,r)$ is defined as follows:
\begin{equation}\relabel{char-fun}
C(M,x)=\sum_{A\subseteq E}(-1)^{|A|}x^{r(M)-r(A)}.
\end{equation}

\begin{exa} 
Let $U_{k,n}$ be the uniform matroid, where $k\le n$, 
i.e.,  the matroid $M=(E,r)$ with $|E|=n$ and 
$r(A)=|A|$ if $|A|\le k$, and $r(A)=k$ otherwise.
By definition, it can be shown that 
$$
C(U_{1,n},x)=
x+\Bra{-{n\choose 1}
	+{n\choose 2}-\cdots+(-1)^n
{n\choose 2}}
=x-1;
$$
$$
C(U_{2,n},x)
=x^2-nx+\Bra{{n\choose 2}-{n\choose 3}+\cdots+
(-1)^n{n\choose n}}
=x^2-nx+(n-1).
$$
Generally, 
for any $1\le k\le n$,
$$
C(U_{k,n},x)=\sum_{i=0}^{k-1} (-1)^{i}{n\choose i} x^{k-i}
+\sum_{i=k}^n (-1)^{i}{n\choose i}.
$$
\end{exa}


\subsection{Computation}

The following properties 
can be derived from the
definition of $C(M,x)$.
Therefore, for any matroid $M$, 
$C(M,x)$ can be determined 
recursively 
by these properties:
\begin{enumerate}[(a)]
	
	\item 
	$C(U_{n,n},x) = (x-1)^n$;
	
	\item if M has a loop, then $C(M,x) = 0$;
	
	\item if $e$ is a coloop of $M$, 
	then $C(M,x) = (x-1)C(M\backslash e,x)$;
	
	\item 
	for matroids 
	$M_1=(E_1,\seti_1)$ and 
	$M_2=(E_2, \seti_2)$
	with $E_1\cap E_2=\emptyset$, 
	the union of $M_1$ and $M_2$, 
	denoted by $M_1\vee M_2$, 
	is the matroid 
	$(E_1\cup E_2, \seti)$,
	where 
	$\seti=\{I_1\cup I_2: I_1\in \seti_1,  I_2\in \seti_2\}$.
	Then, 
	$$
	C(M_1\vee M_2,x)
	=C(M_1,x)C(M_2,x);
	$$
	
	\item if $e$ is not a loop or coloop of $M$, then
	$$
	C(M,x)=C(M\backslash e,x)-C(M/e,x).
	$$
	
\end{enumerate}

For a matroid $M=(E,r)$ and $A\subseteq E$, 
$A$ is called a {\it flat} of  $M$
if $r(A\cup \{e\})>r(A)$ for any $e\in E\setminus A$.
Let $\F(M)$ be the set of flats of $M$ and 
$\mu$ be the m\"obius function $\mu(A,B)$ on flats $A,B$ in $\F(M)$.

Note that $\mu(A,A)=1$ for all $A\in \F(M)$, and 
for each pair of flats $A,B\in \F(M)$
with $A\subset B$:
$$
\mu(A,B)=-\sum\limits_{A\subseteq B'\subset B\atop B'\in \F(M)}\mu(A,B').
$$

\begin{lem} \relabel{char-lem1}
For any flat $F\in \F(M)$, 
$$
\sum_{A\subseteq F\atop r(A)=r(F)} (-1)^{|A|}
=
\left \{
\begin{array}{ll}
\mu(\emptyset,F),  \quad 
&M \mbox{ is loopless};\\
0,&\mbox{otherwise}.
\end{array}
\right.
$$
\end{lem}

\proof Assume that $M$ is loopless.
Define 
$$
U_F=\sum_{A\subseteq F\atop r(A)=r(F)} (-1)^{|A|}.
$$
It is clear that if $F=\emptyset$, then $U_F=1=\mu(\emptyset,F)$.

Now assume that $r(F)\ge 1$. By induction, for any flat $F'<F$
(i.e., $r(F)<r(F)$), the following holds:
$$
U_{F'}=\mu(\emptyset,F').
$$
Then 
\begin{eqnarray*}
\mu(\emptyset,F)
&=& -\sum_{\emptyset \le F'<F}\mu(\emptyset,F')\\
&=&-\sum_{\emptyset \le F'<F}U_{F'}
\\ &=&  -\sum_{\emptyset \le F'<F}\sum_{A\subseteq F'\atop r(A)=r(F')} (-1)^{|A|}
\\ &=&  
\sum_{A\subseteq F\atop r(A)=r(F)}(-1)^{|A|}
-\sum_{\emptyset \le F'\le F}\sum_{A\subseteq F'\atop r(A)=r(F')} (-1)^{|A|}\\
&=&U_F-\sum_{A\subseteq F}(-1)^{|A|}\\
&=&U_F.
\end{eqnarray*}
If $M$ has a loop $e$, then $e\in F$ and the power set $2^F$ 
is partitioned into $2^{F-\{e\}}$ and $\{A\cup \{e\}: A\in 2^{F-\{e\}}\}$.
Thus
$$
U_F=\sum_{A\subseteq F-\{e\}}(-1)^{|A|}+\sum_{A\subseteq F-\{e\}}(-1)^{|A\cup \{e\}|}
=0.
$$
\proofend

Assume that $M$ is loopless, i.e., $r(A)=0$ implies that $A=\emptyset$.

\begin{pro}\relabel{char-pro1}
If $M$ is loopless, then 
$C(M,x)$ has the following expression:
$$
C(M,x)=\sum_{A\subseteq F(M)}\mu(\emptyset,A)x^{r(M)-r(A)}.
$$
\end{pro}

\proof Observe  that 
\begin{eqnarray*}
C(M,x) &=&\sum_{A\subseteq E}(-1)^{|A|}x^{r(M)-r(A)}
\\ &=&
\sum_{F\in \F(M)}
\sum_{A\subseteq F\atop r(A)=r(F)}(-1)^{|A|}x^{r(M)-r(A)}
\\ &=& \sum_{F\in \F(M)}\mu(\emptyset,F)x^{r(M)-r(F)},
\end{eqnarray*}
where the last equality follows from Lemma~\ref{char-lem1}.
\proofend


\subsection{Identity}

\begin{theo}[Kung 2004 \cite{kung2004}]\relabel{char-kung}
For any matroid $M=(E,r)$,
$$
C(M,x_1x_2)=\sum_{F\in \F(M)} C(M/F,x_1)x_2^{r(M)-r(M|F)}C(M|F,x_2),
$$
where $G|F$ is the restriction of $M$ to $F$ 
and $G/F$ is the contraction of $M$ by $F$.
\end{theo}

\proof 
If $F$ is not a flat, then $M/F$ has loops and thus $C(M/F,x)=0$.
Thus the right-hand side can be changed to 
$$
RHS=\sum_{S\subseteq E} C(M/S,x_1)x_2^{r(M)-r(M|S)}C(M|S,x_2).
$$
By definition, 
\begin{eqnarray*}
C(M/S,x)&=&\sum_{S\subseteq A\subseteq E}(-1)^{|A-S|} x^{r(M/s)-r_{M/S}(A-S)}
=\sum_{S\subseteq A\subseteq E}(-1)^{|A-S|} x^{r(M)-r(A)},
\end{eqnarray*}
as $r(M/S)=r(M)-r(S)$ and $r_{M/S}(A-S)=r(A)-r(S)$, while
$$
C(M|S,x)=\sum_{B\subseteq S}(-1)^{|B|}x^{r(S)-r(B)}.
$$
Thus
\begin{eqnarray*}
RHS
&=&\sum_{S\subseteq E} \sum_{B\subseteq S\subseteq A\subseteq E}
(-1)^{|A|+|B|-|S|} x_1^{r(M)-r(A)}x_2^{r(M)-r(B)}           
\\ &=& 
\sum_{B\subseteq A\subseteq E} \sum_{B\subseteq S\subseteq A}
(-1)^{|A|+|B|-|S|} x_1^{r(M)-r(A)}x_2^{r(M)-r(B)}  
\\ &=&  \sum_{B=A\subseteq E} 
(-1)^{|A|} x_1^{r(M)-r(A)}x_2^{r(M)-r(A)}
=C(M,x_1x_2),
\end{eqnarray*}
where the second last equality follows from the following fact that 
if $A$ and $B$ are fixed with $B\subset A$ , then 
$$
\sum_{B\subseteq S\subseteq A}
(-1)^{|A|+|B|-|S|} x_1^{r(M)-r(A)}x_2^{r(M)-r(B)}
=0.
$$
\proofend

\subsection{Connection to other polynomials}

The following relation 
between $C(M,x)$ and $T_M(x,y)$ follows directly 
from their definitions: 
$$
C(M,x)=(-1)^{r(E)}T_M(1-x,0).
$$

For any graph $G$, 
if $M_G$ and $M^*_G$ are the cycle matroid 
and the cocycle matroid of $G$ respectively, 
then 
$$
C(M_G,x)=x^{-c}\chi(G,x), \quad
C(M^*_G , x) = F(G, x).
$$
where $c$ is the number of components of $G$,
and $F(G,x)$ is the flow polynomial of $G$
(see its definition in page~\pageref{eq1-2}).

Thus $C(G,x)$ is an extension of 
both $\chi(G,x)$ and $F(G,x)$.

Due to the fact that 
$C(M_G,x)=x^{-c}\chi(G,x)$, 
the following conclusion 
follows from Proposition~\ref{char-pro1}.

\begin{cor}
	For any simple graph $G=(V,E)$,
	$$
	\chi(G,xy)=\sum_{E'\subseteq E}\chi(G/E',x)\chi(G|E',y),
	$$
	where $G/E'$ is the graph obtained from $G$ by contracting all edges in $E'$
	and $G|E'$ is the graph with edge set $E'$ and vertex set 
	$V_{E'}=\{u\in V: N_u\cap E'\ne \emptyset\}$.
\end{cor}

\subsection{Zero-free interval}

Oxley~\cite{oxl} showed that if every cocircuit 
of $M$ has size at most $d$,  
then $C(M, x) > 0$ holds for all real numbers $r\ge d$.
Jackson~\cite{jac2} pointed out that 
the idea in Oxley's proof can be applied to 
get a more general result. 

A simple minor of $M$ is a minor which contains no loops or circuits of length two.

\begin{theo}[\cite{jac2}]\relabel{C-th1}
Let $M$ be a matroid. 
If every simple minor of $M$ has a cocircuit of size at most $d$,
then $C(M,x) > 0$ for all real numbers $x\ge d$.
\end{theo}

As $F(G,x)=C(M^*_G,x)$, Theorem~\ref{C-th1} implies that 
for any bridgeless graph $G$,
if every 3-edge-connected minor of $G$ has a circuit of 
length at most $d$, then $F(G,x) > 0$ holds 
for all real numbers $t\ge d$. 
It is not difficult to show that  
every 3-connected graph $G$ of order $n$ 
has a circuit of length at most $2 \log_2 n$.
Thus every bridgeless graph of order $n$ has 
all real flow roots less than $2 \log_2 n$.

\subsection{Open Problems}

We conclude this section by 
the following problems. 
A landmark result by Adiprasito, Huh, and Katz
\cite{AHK} shows that for matroids representable over a field, the coefficients of $C(M,x)$ form a log-concave sequence. 
However, the stronger property of real-rootedness remains open in general.

\begin{problem}
	Determine whether 
	$C(M,x)$  is real-rooted for broader classes of matroids
	$M$, or identify precise structural conditions under which real-rootedness holds.
\end{problem}

\begin{problem}
	Are the coefficients of $C(M,x)$ strongly log-concave for all matroids?
\end{problem}

Understanding the location of zeros of $C(M,x)$ is a fundamental problem analogous to that for chromatic polynomials.

\begin{problem}
	Is there a universal zero-free region for $C(M,x)$  depending only on the rank or other invariants of $M$?
\end{problem}

\begin{problem}
	Find direct combinatorial interpretations (e.g., counting certain classes of subsets, orientations, or activities) for the absolute values of the coefficients of 
	$C(M,x)$ for general matroids.
\end{problem}

\section{Chromatic polynomial $\chi(G,x)$}
\relabel{sect-chro}

\def \gch {$\chi(G,x)$}

\subsection{Definition} 

For any simple graph $G$ 
and any positive integer $k$,
a {\it proper $k$-coloring} of $G$
is defined to be a mapping 
$\theta: V(G)\mapsto \{1,2,\cdots,k\}$ 
such that $\theta(u)\ne \theta(v)$ 
for every pair of adjacent vertices $u$ and $v$ in $G$.

The {\it chromatic polynomial} of $G$,
denoted by $\chi(G,k)$,
is defined to be the number of proper $k$-colorings of $G$.

It will be explained later that \gch\ is indeed 
a polynomial in  $x$. 
This graph-function was introduced by G.D. Birkhoff~\cite{Birk1912}
in 1912, 
whose purpose was to use chromatic polynomials as a 
tool to attack the four-colour problem,   
which is equivalent to the statement that 
$\chi(G,4)>0$  for any loopless planar graph $G$.

\begin{exa}
	\begin{enumerate}[(a)]
		\item $\chi(N_n,k)= \lambda^n$ for the null graph $N_n$ of order $n$; 
		
		\item $\chi(T,k)= k
		(k-1)^{n-1}$ for any tree $T$ of order $n$;  and 
		
		\item $\chi(K_n,k)= k(k-1)\cdots (k-n+1)$.
	\end{enumerate} 
\end{exa}

\subsection{Computation} 

\begin{enumerate}[(1)]
	\item {\bf Deletion-contraction formula}:  
	
	\sgap 
	
	\begin{theo} \label{sec5-th1}
		For any simple graph $G$ and any  edge $uv\in E(G)$, 
		\equ{sec5-e2}
		{
			\chi(G,k)=
			\chi(G \setminus uv,k)
			-\chi(G \cdot uv,k),
		}
		where $G \setminus uv$
		is the graph obtained from $G$ by removing edge $uv$
		and 
		$G\cdot uv$ is
		the graph obtained from 
		$G\setminus uv$ by identifying $u$ and $v$ and replacing multiple edges, 
		if any arise, by single ones.
	\end{theo}
	
	The Deletion-contraction formula for $\chi(G,k)$ has another format: 
	for any two non-adjacent vertices $u$ and $v$ in $G$, 
	\equ{sec5-e3}
	{
		\chi(G,k)=
		\chi(G +uv,k)
		+\chi(G \cdot uv,k).
	}

	Thus,  
	for any graph $G=(V,E)$, 
	$\chi(G,k)$ can be determined by the following rules:
	\begin{equation}\label{sec5-e1}
		\chi(G,k)=
		\left \{
		\begin{array}{ll}
			k^{|V|}, &\mbox{if }
			E=\emptyset;\\
			k(k-1)\cdots(k-n+1), &\mbox{if }G \mbox{ is complete};\\
			\chi(G_1,k)\chi(G_2,k),
			&\mbox{if }G=G_1\cup G_2;\\
			\chi(G+uv,k)
			+\chi(G\cdot uv, k),  \quad
			&\mbox{if }u,v\in V, uv\not\in E
		\end{array}
		\right.
	\end{equation}
	where $G_1\cup G_2$ denote the disjoint union of $G_1$ and $G_2$.

	Note that 
	for any real number $k$, or
	even complex number $k$, 
	$\chi(G,k)$ can be determined 
	by (\ref{sec5-e1}). 
	Thus, 
	the function $\chi(G,k)$ can be treated as a function of 
	real numbers, or even complex numbers. 
	
	\sgap 
	
	\item {\bf Chromatic polynomials of chordal graphs}:
	
	\sgap 
	
	A vertex $u$ in a graph $G$ is called a {\it simplicial vertex} of $G$
	if either $d(u)=0$ or 
	$N(u)$ is a clique of $G$.
	Applying (\ref{sec5-e1}),
	one can prove the following result.
	
	\begin{theo} 
		\label{sec5-th2}
		For any simplicial vertex $u$ in a graph $G$,
		$$
		\chi(G,x )=( x-d(u)) \chi(G-u,x ),
		$$
		where $G-u$ is the graph obtained from $G$ by deleting $u$.
	\end{theo}
	
	
	A cycle $C$ in $G$ is called a {\it chordless cycle} of $G$ if 
	$C$ is of length at least $4$, and for any two vertices 
	$u,v$ on $C$, 
	if $uv\in E(G)$, then $uv$ is an edge on $C$.
	If $G$ does not contain any chordless cycle, then $G$ is called a {\it chordal graph}.

	A {\it perfect elimination ordering} of  $G$ is an ordering $v_1,v_2,\cdots, v_n$ of vertices in $G$  
	such that $v_i$ is a simplicial vertex of 
	the subgraph of $G$ induced by 
	$\{v_1,v_2,\cdots,v_i\}$ 
	for all $i=1,2,\cdots,n$. 
	
	\begin{theo}[Dirac~\cite{Dirac1961}]
		\label{sec5-th3}
		A graph $G$ is chordal if and only if it has a perfect elimination ordering.
	\end{theo} 
	
	For any chordal graph $G$, 
	if $v_1,v_2,\cdots,v_n$ is 
	a perfect elimination ordering
	of $G$,
	applying Theorem~\ref{sec5-th2}
	yields that 
	\equ{sec5-e4-0}
	{
		\chi(G,x)=\prod_{i=1}^n 
		(x-d_{G_i}(v_i)),
	}
	where $G_i$ is the subgraph of $G$ induced by 
	$\{v_1,v_2,\cdots,v_i\}$.
	
	The following conclusion follows directly from Theorem~\ref{sec5-th2}.
	
	\begin{cor}\label{sec5-cor1}
		For any chordal graph $G$ 
		with chromatic number $\chi(G)=r$, we have 
		\equ{sec5-e4}
		{
			\chi(G,x)=x^{n_0}(x-1)^{n_1}
			\cdots (x-r+1)^{n_{r-1}},
		}
		where $n_i$ is a positive integer 
		for each $i=0,1,\cdots,r-1$.
	\end{cor} 
	
	
	\item {\bf Chromatic polynomial
		of the join of two graphs}:
	
	\sgap 
	
	For any two vertex-disjoint graphs $G$ and $H$,
	the {\it join} of them, \label{chap5-page-join}
	denoted by $G+H$,  
	is defined to be the graph 
	with vertex set $V(G)\cup V(H)$ and edge set
	$E(G)\cup E(H)\cup \{xy: x\in V(G),y\in V(H)\}$.
	An expression of $\chi(G+H,\lambda)$ in terms of 
	$\chi(G,\lambda)$ and $\chi(H,\lambda)$
	was provided by Zykov~\cite{zyk1949}.

	\begin{theo}\relabel{chap5-factor4}
		For any two graphs $G$ and $H$, if 
		$$
		\chi(G,x)
		=\sum_{i=1}^{n_1}
		a_i(x)_i\quad
		\mbox{and}\quad
		\chi(H,x)
		=\sum_{j=1}^{n_2}b_j
		(x)_j,
		$$
		then $\chi(G+H,x)=
		\sum\limits_{i=1}^{n_1}\sum\limits_{j=1}^{n_2}a_ib_j
		(x)_{i+j}$, 
		where $(x)_i=x(x-1)\cdots (x-i+1)$.
	\end{theo}
	
	For any integers $n$ and $k$ with $0\le k\le n$,  the
	{\it Stirling number of the second kind}, 
	denoted by $S(n,k)$,
	\index{Stirling numbers of the second kind} 
	\red{counts the number of 
		ways to partition a set of $n$ 
		objects into $k$ non-empty subsets.}
	The function $S(n,k)$ satisfies the following 
	recursive expression: for any $2\le k\le n-1$,
	$$
	S(n,k)=S(n-1,k-1)+kS(n-1,k).
	$$
	
	The definitions of $S(n,k)$ and $\chi(N_r,\lambda)$ 
	imply that 
	$$
	\chi(N_r,x)
	=x^r=
	\sum_{i=1}^r S(r,i)(x)_i.
	$$
	For any complete $k$-partite graph $G$,
	Theorem~\ref{chap5-factor4} gives 
	an explicit expression for $\chi(G,x)$.
	
	\begin{cor}\relabel{chap5-factor5}
		For the complete $k$-partite graph 
		$K(p_1,p_2,\cdots, p_k)$,
		\begin{equation}
			\relabel{chap5-factor5-eq1}
			\chi(K(p_1,p_2,\cdots, p_k),x)
			=\sum_{r_k=1}^{p_k}\cdots \sum_{r_1=1}^{p_1}
			\left (\prod_{i=1}^k S(p_i,r_i)\right )
			(x)_{r_1+\cdots+r_k}.
		\end{equation}
	\end{cor}
	
\end{enumerate} 

\subsection{Expansion} 


Applying (\ref{sec5-e1}), the following can be established. 

\begin{theo}\label{sec5-th6}
	Let $G$ be a simple graph of
	order $n$. Then
	$\chi(G,x)$ is a polynomial in $x$ of degree $n$, as expressed below:
	\equ{sec5-e7}
	{
		\chi(G,x)=\sum_{i=0}^n a_ix^n; 
	}
	where $a_n=1$, 
	$a_i=0$ for each $i$ with $i<c(G)$ (in particular, $a_1=0$),
	each $a_i$ is an integer
	and $a_i$ has the sign 
	$(-1)^{n-i}$.
\end{theo}

In 1932, Whitney established the following result which 
interprets the coefficients $a_i$ 
of $x^i$ in the expansion of $\chi(G,x)$. 

\begin{theo}[Whitney~\cite{whi1932}]
	\label{sec5-th4}
	For any simple graph 
	$G=(V,E)$, 
	\begin{equation}
		\label{sec5-e5}
		\chi(G,x)=
		\sum_{A\subseteq E}
		(-1)^{|A|}x^{c(A)},
	\end{equation}
	where $c(A)$ is the number of components of $(V,A)$ (i.e.,  
	the spanning subgraph of $G$ with edge set $A$).
\end{theo}

For any integer $r$ with $1\le r\le |V|$,  let  
${\cal E}_r$ denote the family of $A\subseteq E$ 
with $c(A)=r$.  Then, 
(\ref{sec5-e5}) implies that 
\begin{equation}
	\label{sec5-e8}
	\chi(G,x)=
	\sum_{r=1}^{|V|}
	\Bra{\sum_{A\in {\cal E}_r}
		(-1)^{|A|}} 	x^{r}.
\end{equation}

Let $G$ be a graph with $n$ vertices 
and $m$ edges, and let  
$\beta:E(G)\mapsto \{1,2,\cdots,m\}$ be 
a bijection. 
For any cycle $C$ in $G$, let
$e$ be the edge in $E(G)$ such that $\beta(e)\le \beta(e')$ for every $e'$ in $E(C)$. 
We call the path $C-e$ a 
{\it broken-cycle}
in $G$ with respect to $\beta$.

\begin{theo}[Whitney's Broken-cycle Theorem]
	\label{sec5-th7}
	Let $G$ be a graph with 
	$n$ vertices and $m$ edges, and let 
	$
	\beta:E(G)\mapsto \{1,2,\cdots,m\}
	$ 
	be any bijection.
	Then 
	\begin{equation}
		\label{sec5-e9}
		\chi(G,x)=\sum_{i=1}^n (-1)^{n-i}h_i(G) x^i,
	\end{equation}
	where $h_i(G)$ is the number of spanning subgraphs of 
	$G$ that have exactly $n-i$ edges and contain
	no broken cycles with respect to $\beta$.
\end{theo}

For any graph $H$, let $n_G(H)$
(resp. $i_G(H)$) denote the 
number of subgraphs (resp.
induced subgraphs) of $G$ which are isomorphic to $H$.
If the graph $G$ under consideration is clear, 
$n_G(H)$ and $i_G(H)$
will be replaced by $n(H)$ 
and $i(H)$.
By Theorem~\ref{sec5-th7}, 
we obtain the following 
interpretations for some
coefficients of $\chi(G,x)$
(see \cite{don2005}). 

\begin{cor}\label{sec5-cor2}
	Let $G$ be a simple graph of order $n$, size $m$ and girth $g$. Then 
	\begin{equation}
		\label{sec5-e10}
		\chi(G,x)=
		\sum_{i=1}^{n} (-1)^{n-i}h_i x^i,
	\end{equation}
	where 
	\begin{enumerate}[(a)]
		\item for $0\le i\le g-2$,
		$h_{n-i}={m\choose i}$
		(in particular,  $h_n=1$ and 
		$h_{n-1}=m$);
		\item $h_{g-1}={m\choose g-1}-n(C_g)$
		(in particular,  
		$h_{n-2}={m\choose 2}-n(C_3)$),
		where $n(C_g)$ is the number of subgraphs in $G$ which are isomorphic to 
		a cycle of length $g$; and 
		
		\item $h_{n-3}={m\choose 3}-(m-2)n(K_3)-i(C_4)+2n(K_4)$.
		
	\end{enumerate} 
\end{cor}

\subsection{Expressions in terms of acyclic orientations}

In this subsection, let $G$ be a simple graph of order $p$.
In 1973, Stanley established 
a relation between $\chi(G,k)$ 
and acyclic orientations of $G$,
as stated below.

\begin{pro}[Stanley \cite{sta1973a}]\relabel{pro1-1}
	For any non-negative integer $k$, \red{$\chi(G,k)$} is equal to the number of ordered pairs 
	$(\theta,\O)$, where $\theta$ is a mapping
	$\theta: V \mapsto  \{1, 2,\cdots, k\}$ and $\O$ is an
	orientation of $G$, subject to the two conditions below:
	\begin{enumerate}[(a)]
		\item the orientation $\O$ is acyclic; and 
		
		\item for every arc 
		$u \rightarrow v$ 
		in the orientation $\O$, 
		$\theta(u)<\theta(v)$ holds.
	\end{enumerate} 
\end{pro}

\proof Define a mapping $\psi$ with $\psi(f)=(\theta,\O)$ from the set of 
proper $k$-colourings $f$ of $G$ 
to the set of ordered pairs $(\theta,\O)$, 
where  
$\theta=f$ and $\O$ is the orientation of $G$ such that 
$u\rightarrow v$ is an arc whenever $uv\in E(G)$ and $f(u)<f(v)$.

Clearly, for any given $f$, the ordered pair $(\theta,\O)$ defined above 
satisfies conditions (a) and (b).
It is obvious that $\psi$ is a bijection and thus the result holds.
\proofend

Let {\bf $\tchi(G, k)$} denote 
the number of ordered pairs 
$(\theta,\O)$, where $\theta$ is any mapping 
$\theta: V \mapsto  \{1, 2,\cdots, k\}$ and $\O$ is an
orientation of $G$, subject to the two conditions below:
\begin{enumerate}[(a)]
	\item the orientation $\O$ is acyclic; and 
	
	\item for every arc
	 $u \rightarrow v$ 
	in the orientation $\O$, 
	$\theta(u)\le \theta(v)$ holds.
\end{enumerate} 

The relationship between  $\chi(G, x)$ and 
$\tchi(G, x)$ is somewhat analogous to the relationship between combinations of $n$ things taken
$k$ at a time without repetition, enumerated by
${n\choose k}$, and with repetition, enumerated by 
${n+k-1\choose k}=(-1)^k {-n\choose k}$. 
(\blue{Note that 
	${n+k-1\choose k}$ is the number of non-negative integer 
	solutions of $x_1+x_2+\cdots+x_n=k$.})

\begin{theo}[Stanley \cite{sta1973a}]\relabel{theo1-1}
	For any graph $G$ of order $p$
	and  
	non-negative integer $x$,
	$$
	\tchi(G, x)=(-1)^p \chi(G, -x), \quad i.e.,\quad 
	\chi(G, x)=(-1)^p \tchi(G, -x).
	$$
\end{theo}

\proof It suffices to show that 
\begin{enumerate}[(a)]
	\item $\tchi(N_1, x)=x$;
	
	\item  $\tchi(G\cup H, x)=
	\tchi(G, x)\tchi(H, x)$;
	and 
	
	\item $\tchi(G, x)=\tchi(G\backslash e, x)
	+\tchi(G/e, x)$ holds for any edge $e$.
	\proofend
\end{enumerate}

Theorem~\ref{theo1-1} provides a combinatorial interpretation of the positive integer $(-1)^p\chi(G,-k)$, 
where $k$ is a positive integer. 
In particular, when $k=1$, 
every orientation of $G$ is automatically compatible with every map $\theta : V\mapsto  \{1\}$.

\begin{cor}\relabel{sec5-cor3} 
	If $G$ is a graph with $p$ vertices, then 
	$(-1)^p\chi(G,-1)$ is equal to the number of 
	acyclic orientations of $G$.
\end{cor}

Let $\seta(G)$ be the set of acyclic orientations of $G$.
For any $D\in \seta(G)$, 
let 
$\Sorder(D,k)$  (resp. $\order(G,k)$)
be the number of 
{\it strictly order-preserved mappings} 
(resp. {\it order-preserved mappings})
$\theta: V(G)\mapsto \{1,2,\cdots,k\}$
with respect to $D$,
i.e., $\theta(u)<\theta(v)$ 
(resp. $\theta(u)\le \theta(v)$) 
holds for every arc
$u\rightarrow v$ in $D$.

By Proposition~\ref{pro1-1},
the following equality holds:
\begin{equation}\relabel{eq1-3}
	\chi(G,k)=\sum_{D\in A(G)}\Sorder(D,k).
\end{equation}

For any $D\in \seta(G)$, 
a permutation of $\pi=(i_1,i_2,\cdots, i_p)$ of $1,2,\cdots,p$ is said to be 
{\it order preserved}
if $i_s$ appears before $i_t$ (i.e., $s<t$) in $\pi$ 
for every arc 
$i_s\rightarrow i_t$  in $D$.
Let $\OP(D)$ be the set of all
order preserved permutations 
$\pi=(i_1,i_2,\cdots, i_p)$ of $1,2,\cdots,p$.
For any $\pi=(i_1,i_2,\cdots, i_p)\in \OP(D)$, let 
$\rho(\pi)= 
|\{1\le j\le p-1: i_j<i_{j+1}\}|$.
Obviously, $0\le \rho(\pi)\le p-1$.

Stanley \cite{sta1970} showed that for any $D\in \seta(G)$, 
\equ{sec5-e11}
{
	\Sorder(D,k)=\sum_{\pi \in \OP(D)} {k+p-1-\rho(\pi)\choose p}.
}

By (\ref{eq1-3}) and 
(\ref{sec5-e11}), the following result is obtained.

\begin{theo}[Stanley \cite{sta1970}]
	\relabel{theo1-2}
	$$
	\chi(G,x)=\sum_{D\in \seta(G)}
	\sum_{\pi \in \OP(D)} {x+p-1-\rho(\pi)\choose p}.
	$$
\end{theo}


Similarly, Stanley \cite{sta1970} also established the following conclusion.

\begin{theo}[Stanley \cite{sta1970}]
	\relabel{theo1-3}
	$$
	\tchi(G,x)
	=\sum_{D\in \seta(G)} \order(D,x)
	=\sum_{D\in \seta(G)}
	\sum_{\pi \in \OP(D)} {x+\rho(\pi)\choose p}.
	$$
\end{theo}

Let $G=(V,E)$ be a simple graph with $V=\{1,2,\cdots,p\}$.
For any permutation  $\pi=(v_1,v_2,\cdots,v_p)$ of 
$\{1,2,\cdots,p\}$,
let $\delta_G(\pi)$
(or simply $\delta(\pi)$) be the number of $i$'s, where $1\le i\le p-1$, such that either $v_i<v_{i+1}$ or $v_iv_{i+1}\in E$.
Let $\W(G)$ be the set of subsets $\{a,b,c\}$ of $V$,
where $a<b<c$, 
such that $G[\{a,b,c\}]$ contains 
$ac$ as its only edge.


\begin{theo}[Dong \cite{don2020}]\label{dong-th3}
	For any simple graph $G=(V,E)$
	with $|V|=p$, 
	$\W(G)=\emptyset$  if and only if  
	$(-1)^p\chi(G,-x)=\sum\limits_{\pi} 
	{x+\delta(\pi)\choose p}$, 
	where the sum runs over all $p!$ permutations  of $V$.
\end{theo}

\def \NW {\mathscr{NW}}

\subsection{Mean color number}

\def \meanc {{\mu}}

Let $G=(V,E)$ be a graph of order $n$.
The {\it mean color number} of $G$, denoted by $\meanc(G)$, is defined as the average number of colors used across all possible proper 
$n$-colorings of $G$
(see \cite{Bartels}).
For any $1\le k\le n$,
let $\alpha_k(G)$ denote the number
of ways to partition $V$ into exactly
$k$ nonempty independent sets.
By definition,
\equ{eq3-0}
{\chi(G,x)
=\sum\limits_{k=1}^n \alpha_k(G) (x)_k,
}
where $(x)_k=x(x-1)\cdots (x-k+1)$, 
and 
\Eqn{eq3-1}
{
	\meanc(G) &=&
	\frac{ 
		\sum_{k=1}^n \alpha_k(G)\,
		k\, (n)_k
	}
	{
		\sum_{k=1}^n \alpha_k(G) (n)_k
	}\nonumber \\
	&=&
	\frac{ 
		\sum_{k=1}^n 
		\alpha_k(G) 
		(n- (n-k))
		(n)_k
	}
	{
	\chi(G,n)
	}\nonumber  \\
&=&\frac{n\chi(G,n)-	n\sum_{k=1}^n 
	\alpha_k(G) 
	(n-1)_{k}}
{\chi(G,n)}
\nonumber 
\\
&=&n
\Bra{1-\frac{\chi(G,n-1)}{\chi(G,n)}
},
}
as shown in 
 \cite{Bartels}.
Applying (\ref{eq3-1}),
an explicit expression of  $\meanc(G)$ can be obtained 
when $G$ is 
 the complete graphs
$K_n$ or the null graph $N_n$:
\Enum{[(i)]
\item $\meanc(K_n)=n$; and 
\item $\meanc(N_n)=
n\Bra{1-\frac{(n-1)^n}{n^{n}}}$.
}

Clearly, for any graph $G$ of order $n$,
\equ{eq3-3}
{
\meanc(G)
\le n
=\meanc(K_n).
}
Is it true in general that 
a graph with more edges tend to have a larger mean color number?
Bartels and Welsh~\cite{Bartels}
conjectured that 
$\meanc(H)\le \meanc(G)$
holds
for any spanning subgraph $H$
of $G$.
However, Mosca~\cite{Mosca}  found counterexamples to this conjecture.

Another conjecture proposed by Bartels and Welsh~\cite{Bartels} is the so-called ``Shameful Conjecture", stated below,
which has been proven.

\begin{con}[\cite{Bartels}]
	\label{con3-2}
	For any graph $G$ of order $n\ge 2$, 
\equ{eq3-4}
{
\meanc(G)\ge 
n\Bra{1-\frac{(n-1)^n}{n^{n}}}
=\meanc(N_n).
}
\end{con}

By (\ref{eq3-1}) and 
(\ref{eq3-4}), 
$\meanc(G)\ge 
n\Bra{1-\frac{(n-1)^n}{n^{n}}}$
if and only if 
\equ{eq3-5}
{
\frac{\chi(G,n)}{\chi(G,n-1)}
\ge \frac{n^n}{(n-1)^{n}}.
}
Note that $\frac{n^n}{(n-1)^{n}}$ 
is strictly decreasing 
for $n\ge 2$
and tends to $e$ when $n$ goes to infinity,
where $e$ $(=2.7182818\cdots)$ is the base of natural
logarithms.
Thus,  Conjecture~\ref{con3-2} 
implies the
following inequality 
\equ{eq3-6}
{
\frac{\chi(G,n)}{\chi(G,n-1)}\ge e.
}
Seymour~\cite{Seymour} obtained the following result, which is very close to 
(\ref{eq3-6}):
\equ{eq3-7}
{
	\frac{
\chi(G,n)}
{\chi(G,n-1)}
\ge \frac{685}{252}
\, (= 2.7182539\ldots ).
}
Dong~\cite{Dong2000} proved
 that 
\equ{eq3-80}
{
\frac{\chi(G,x)}
{\chi(G,x-1)}\ge 
\frac{x^n}{(x-1)^n}
	}
for all real numbers 
$x\ge n$. 
Clearly, (\ref{eq3-5}) 
is a special case 
of (\ref{eq3-80}).

Bartels and Welsh~\cite{Bartels} also 
proposed the following conjecture which is much 
stronger than Conjecture~\ref{con3-2}.

\begin{con}
	[\cite{Bartels}]
	\label{con3-1}
	If $G$ is not an empty graph,
	there exists an edge $uv$ 
	in 
	$G$ such that 
	$\meanc(G-uv)\le \meanc(G)$.
\end{con}

As far as I know, Conjecture~\ref{con3-1} remains open,
although it holds for some families of graphs.

\subsection{Known results}

In this subsection, we present some established results on chromatic polynomials, especially concerning the locations of chromatic roots. Further results are available in \cite{don2022, don2005}.

\begin{theo}\label{sec5-th13}
	\begin{enumerate}[(a)]
		\item $\chi(G,x)\ne 0$
		for any graph $G$ and 
		$x\in (-\infty,1)\setminus \{0\}$;
		
		\item{}  
		(Jackson \cite{jac1993})
		$\chi(G,x)$ is non-zero 
		for any graph $G$ and any 
		$x\in (1,32/27]$; 
		
		\item{} (Thomassen \cite{tho1997}) for any real numbers $x_1$ and $x_2$ with $32/27\le x_1<x_2$, there exists a graph $H$ 
		a number $x\in (x_1,x_2)$
		such that $\chi(H,x_0)=0$;
		
		\item{} (Sokal \cite{sok2001})
		for any graph $G$,
		$\chi(G,z)\ne 0$
		for any complex number $z$
		with $\|z\|\ge  C\Delta(G)$,
		where  $C\le 7.963907$,
		and this result was strengthened 
		to 
		that  $C\le 6.907 \cdots$
		by Fern\'andez and  Procacci \cite{fer2008};
		
		\item{} (Sokal \cite{sok2004})
		roots of chromatic polynomials of graphs 
		are dense in the whole complex plane;

		\item{} (Royle \cite{roy2008}) for any $\epsilon>0$, there exists 
		a plane triangulation $X$ 
		such that $\chi(X,x)$ has a 
		real root in the interval $(4-\epsilon,4)$;
		
		\item{} (Tutte \cite{tutte1970})
		let $\varphi$ denote
		the Golden Ratio $\frac{1+\sqrt 5}2$.
		Then $\chi(T, \varphi+1)\ne 0$
		and $\chi(T, \varphi+2)\ne 0$
		for every planar triangulation 
		$T$.
		This conclusion has been extended to 
		all planar graphs 
		by Perrett and Thomassen 
		\cite{per2018}; and 
		
		\item{} (Birkhoff and Lewis \cite{Birk1946}) $\chi(G,x)>0$ 
		for any loopless planar graph $G$ and  real number
		$x\in [5,\infty)$.
		
	\end{enumerate} 
	
\end{theo}

\begin{theo}\label{sec5-th14}
	Let $G$ be a simple graph of order $n$, size $m$, 
	maximum degree $\Delta$
	and chromatic number $r$.
	Assume that $\chi(G,x)
	=\sum\limits_{i=1}^n(-1)^{n-i} h_ix^i$.  Then
	\begin{enumerate}[(a)]	
		\item{} (Stanley \cite{sta1973a})
		$|\chi(G,-1)|$ equals
		to the number of acyclic orientations of $G$; 
		
		\item{} (Greene and Zaslavsky \cite{gre1983})
		for any $v\in V(G)$, 
		$h_1$ equals to the number of acyclic orientations 
		of $G$ with $v$ as its unique source;

		\item{} (Huh \cite{huh2012})
		$h_{i-1}h_{i+1}\le h_i^2$
		for all $i=2,3,\cdots,n-1$,
		and this result has been extended to 
		$C(M,x)$ for all 
		representable matroids $M$
		by Huh and Katz
		\cite{huh2012b};
		
		\item (Dong and Zhang \cite{dong2023})
		$\chi(G,k)=P_{\ell }(G,k)$
		holds for all integer $k$ with $k\ge m-1$, 
		where $P_{\ell }(G,k)$ is
		the list color function of $G$, which is defined to be the minimum number of proper  $L$-colorings 
		over all $k$-assignments $L$; 
		and 
		
		\item{} (Brown and Erey~\cite{Brown2015}) 
		$	\chi(G,x)\le (x)_r (x-1)^{\Delta-(r-1)}x^{n-1-\Delta}$	holds 
		for any integer $x\ge 0$.
		
	\end{enumerate} 
	
\end{theo}

\subsection{Open problems}

Here we list some well-known problems on chromatic polynomials, and more 
open problems can be 
found in \cite{don2022, don2005}.

In 1946, Birkhoff and Lewis~\cite{Birk1946} 
proved that $\chi(G,x)> 0$ 
for all real numbers $x\in [5,\infty)$.  They also conjectured that this conclusion holds for all real $x\in [4,5)$.

\begin{con}[Birkhoff and Lewis~\cite{Birk1946}]
	\label{sec5-con1}
	For any planar and loopless graph $G$, 
	$\chi(G,x)>0$ holds for all real numbers $x\in [4,5)$.
\end{con} 

So far the only known fact on 
Conjecture~\ref{sec5-con1}
is the case $x=4$ (i.e., the four-color Theorem).

For any integer $k\ge 1$,
let $\Omega_k$ be the set of roots of 
monic polynomials of degree $k$ 
whose coefficients are integers.
Two conjectures related to $\Omega_k$,
called 
{\it the $\alpha+n$ conjecture} 
\index{the $\alpha+n$ conjecture}
and {\it the $\alpha n$ conjecture}, 
\index{the $\alpha n$ conjecture}
were proposed during 
the programme on ``Combinatorics and Statistical
Mechanics" 
at the Newton Institute of Cambridge in 2008
(\cite{cam2017}).

\begin{con}[\cite{cam2017}]
	\label{sec5-con2}
	\begin{enumerate}[(a)]
		\item{} 
		(The $\alpha+n$ conjecture)
		For every $\alpha\in \Omega_k$, where $k\ge 1$,
		there exists an integer $n$ 
		such that $\alpha+n$ is a chromatic root;
		
		\item{} (The $\alpha n$ conjecture)
		Let $\alpha$ be a chromatic root. 
		Then $n\alpha$ is also a chromatic root for 
		every positive integer $n$.
	\end{enumerate} 
\end{con}

In 2001, Sokal~\cite{sok2001} 
proved that for any graph $G$, 
all complex roots 
of $\chi(G,x)$ are within a disc 
$\{z: \|z\|\le C\Delta(G)\}$, where 
$C\le 7.963907$.
A stronger conjecture
on complex chromatic roots
was proposed by 
Sokal~\cite{sok2000} in 2000.


\begin{con}[Sokal~\cite{sok2000}] 	\label{sec5-con3}
	For any graph $G$ with maximum degree $\Delta$, 
	if $\chi(G,z)=0$, then 
	$\|z\|\le 2\Delta$.
\end{con} 

For real roots of $\chi(G,x)=0$,
he gave the following conjecture. 

\begin{con}[Sokal~\cite{sok2005}] \label{sec5-con4}
	For any graph $G$ with maximum degree $\Delta$,
	if $\chi(G,x)=0$ for real number $x$, then 
	$x\le \Delta$.
\end{con}

The ``maximaxflow" 
denoted by $\Lambda (G)$, is defined below:
$$
\Lambda (G)=\max_{u\ne v\atop u,v\in V(G)} \lambda(u,v),
$$
where 
$\lambda(u,v)$ is the maximum number of edge-disjoint 
$u-v$ paths, which is equal to the 
minimum number of edges separating $u$ from $v$.
Obviously, $\lambda(u,v)\le 
\min\{d_G(u),d_G(v)\}$
and $\Lambda(G)\le \Delta(G)$.
Sokal~\cite{sok2005} further conjectured the following.

\begin{con} [Sokal \cite{sok2005}]
	\label{sec5-con5}
	For any graph $G$, if 
	$\chi(G,x)=0$ for real number $x$, then 
	$x\le \Lambda (G)$.
\end{con}

The following question of Beraha~\cite{ber1975} about real chromatic roots of planar
graphs remains unsolved  for
$n = 8, 9$ and $n \ge 11$.

\begin{prob}[Beraha~\cite{ber1975}] \label{sec5-prob3}
	Let $B_n = 2 + 2 cos(2\pi/n)$ 
	for $n\in \Z$. Is it true that for every 
	$\epsilon>0$, there exists a plane triangulation $G$ such
	that $\chi(G,x)$
	has a root in $(B_n -\epsilon, 
	B_n + \epsilon)$?
\end{prob} 

Thomassen~\cite{tho1996} conjectured that 
Hamiltonian graphs have no chromatic roots in $(1,2)$.

\begin{con}
	[Thomassen \cite{tho1996}]
	\label{sec5-con6}
	If $G$ is a Hamiltonian graph,
	then $\chi(G,x)\ne 0$ for all real $x\in (1,2)$.
\end{con}

A connected graph $G$ is said to be {\it $\alpha$-tough} 
if $G-S$ has at most $|S|$ 
components 
for any non-empty 
independent set $S$ of $G$.
Obviously, any Hamiltonian graph is $\alpha$-tough.
Conjecture~\ref{sec5-con6} has been strengthened 
by replacing a Hamiltonian graph $G$ to a $\alpha$-tough graph $G$ by Dong and Koh~\cite{don2006}.

A graph $G$ is said to be {\it $\chi$-unique} if for any graph $H$, $\chi(H,x)=\chi(G,x)$ 
always 
implies that $H\cong G$.
Bollob\'{a}s, Pebody and Riordan
~\cite{bollo2000} proposed the following conjecture
on the trend of proportion of $\chi$-unique graphs.

\begin{con}
	[\cite{bollo2000}]
	$F_{\chi}(n) \rightarrow  1$ 
	as $n \rightarrow  \infty$, where $F_{\chi}(n)$ is the proportion of $\chi$-unique graphs
	of order $n$. 
\end{con}

For positive integers $r_0,r_1,\cdots,r_k$,
where $k\ge 0$, 
let ${\cal I}(r_0,r_1,\cdots,r_k)$ 
denote the family of graphs $G$
with 
$\chi(G,x)=x^{r_0}(x-1)^{r_1}\cdots 
(x-k)^{r_k}$.

\begin{prob}
	\label{sec5-prob2} 
	Given $k$ and $ r_0,r_1,\cdots,r_k\in 
	\Z$, characterize the graphs 
	in the family ${\cal I}(r_0,r_1,\cdots,r_k)$. 
	In particular, determine if the family ${\cal I}(r_0,r_1,\cdots,r_k)$
	consists of chordal graphs
	only.
\end{prob}

For any simple graphs $G$ and $H$, 
the {\it tensor product of $G$ and $H$}, denoted by $G\times H$,  
is the graph such that its 
vertex set is 
$
\{(u,v): u\in V(G), v\in V(H)\}
$
and two vertices $(u_1,v_1)$ and $(u_2,v_2)$ are adjacent in $G\times H$ if and only if 
\begin{itemize} 
	\item	$u_1$ and $u_2$ 
	are adjacent in $G$; and 
	\item 
	$v_1$ and $v_2$ 
	are adjacent in $H$.
\end{itemize} 

For example, $C_5\times K_2$ 
is the graph shown in 
Figure~\ref{jimei-f5}.
It is not difficult to verify that 
for any graph $G$, 
$G\times K_2$ is a bipartite graph of order $2|V(G)|$
and size $2|E(G)|$.

\begin{figure}[h!]
	\centering

\begin{tikzpicture}[
	line width=0.85pt,
	line cap=round,
	line join=round,
	vertex/.style={circle,draw=black,fill=black,
		inner sep=0pt,minimum size=6.8pt},
	vertex label/.style={font=\Large,inner sep=1pt},
	pair label/.style={font=\large,inner sep=1pt},
	graph label/.style={font=\LARGE}
	]
	
	\begin{scope}
		\coordinate (u1) at ( 0,     1.50);
		\coordinate (u2) at ( 1.14,  0.49);
		\coordinate (u3) at ( 0.53, -0.82);
		\coordinate (u4) at (-0.74, -0.80);
		\coordinate (u5) at (-1.30,  0.53);
		
		\draw (u1) -- (u2) -- (u3) -- (u4) -- (u5) -- cycle;
		\foreach \v in {u1,u2,u3,u4,u5}
		\node[vertex] at (\v) {};
		
		\node[vertex label,anchor=south] at (0,1.65) {$u_1$};
		\node[vertex label,anchor=west] at (1.30,0.49) {$u_2$};
		\node[vertex label,anchor=west] at (0.69,-0.91) {$u_3$};
		\node[vertex label,anchor=east] at (-0.92,-0.85) {$u_4$};
		\node[vertex label,anchor=east] at (-1.47,0.57) {$u_5$};
		
		\node[graph label] at (0,-2.02) {$C_5$};
	\end{scope}
	
	\begin{scope}[xshift=2.825cm]
		\coordinate (v1) at (0,0.81);
		\coordinate (v2) at (0,-0.71);
		
		\draw (v1) -- (v2);
		\node[vertex] at (v1) {};
		\node[vertex] at (v2) {};
		\node[vertex label,anchor=west] at (0.16,0.81) {$v_1$};
		\node[vertex label,anchor=west] at (0.16,-0.71) {$v_2$};
		
		\node[graph label] at (0.2,-2.02) {$K_2$};
	\end{scope}
	
	\begin{scope}[xshift=7.725cm]
		\coordinate (p11) at ( 0,       1.6625);
		\coordinate (p12) at ( 0,       1.0875);
		\coordinate (p21) at ( 1.6500,  0.4250);
		\coordinate (p22) at ( 1.2125,  0.2500);
		\coordinate (p31) at ( 0.9875, -1.4000);
		\coordinate (p32) at ( 0.6375, -0.9000);
		\coordinate (p41) at (-0.9500, -1.4000);
		\coordinate (p42) at (-0.7875, -0.8875);
		\coordinate (p51) at (-1.7375,  0.4250);
		\coordinate (p52) at (-1.2375,  0.3250);
		
		\draw (p11) -- (p22) -- (p31) -- (p42) -- (p51)
		-- (p12) -- (p21) -- (p32) -- (p41) -- (p52) -- cycle;
		
		\foreach \v in {p11,p12,p21,p22,p31,p32,p41,p42,p51,p52}
		\node[vertex] at (\v) {};
		
		\node[pair label,anchor=south] at (0.18,1.83) {$(u_1,v_1)$};
		\node[pair label,anchor=north] at (-0.02,0.86) {$(u_1,v_2)$};
		\node[pair label,anchor=west] at (1.82,0.46) {$(u_2,v_1)$};
		\node[pair label,anchor=east] at (1.02,0.15) {$(u_2,v_2)$};
		
		\node[graph label] at (0.3,-2.08) {$C_5\times K_2$};
	\end{scope}
	
\end{tikzpicture}
	
	\caption{The graph $C_5\times K_2$
	\relabel{jimei-f5}
}
\end{figure}

In 2011,  Zhao~\cite{zhao2011}
proposed the following conjecture which has been settled 
for all $k\ge (2n)^{2n+2}$
by himself, where $n=|V(G)|$. 

\begin{con}[\cite{zhao2011}]
	\label{zhao2011}
	For any simple graph $G$,
$
\chi(G\times K_2,k)
\ge \chi(G,k)^2
$ 
for all $k\in \Z$.
\end{con}

The unimodality conjecture for chromatic polynomials, originally proposed by Read~\cite{read1968}, was resolved by Huh~\cite{huh2012} (see Theorem~\ref{sec5-th14} (c)).
Nevertheless, this well-known conjecture deserves a purely combinatorial proof, 
as its formulation is inherently combinatorial. For this reason, we restate the problem here.

\begin{prob}
	\label{uni-prob}
Is there a combinatorial proof for the unimodal conjecture
for chromatic polynomials?
More precisely,
for any graph $G$ of order $n$, 
if 
$$
\chi(G,x)
=\sum\limits_{1\le i\le n}(-1)^{n-i}h_ix^i, 
$$
does there exist a combinatorial proof of this log-concavity property of the sequence 
$\{h_i\}_{1\le i\le n}$?
\end{prob}

We end this section by the following problem,
proposed by Wilf and Kinial,
which was 
recently surveyed by Lazebnik~\cite{laz2019}.

\begin{prob}[\cite{laz2019}]
	\label{sec5-prob1}
	Let $m, n,k\in \Z$, where $m\le {n\choose 2}$.
	What is the maximum value of $\chi(G,x)$, where $x\le k$, 
	over all graphs $G$ 
	of order $n$ and size $m$?
\end{prob}

\section{Flow polynomial $F(G,x)$}
\relabel{sect-flow}

\subsection{Flow problems}

 Let $D$ be any orientation of a graph $G$
and $\Gamma$ be any Abelian group.
Let $f$ be a mapping 
$f:A(D)\mapsto \Gamma$, where $A(D)$ 
is the set of arcs in $D$.
$f$ is called a {\it $\Gamma$-flow} of $D$ 
if the following equality holds 
at every vertex $u$ of $D$:
\equ{sec4-e1}
{
\sum_{a\in A^+(u)}f(a)=\sum_{a\in A^-(u)}f(a),
}
where $A^+(u)$ (resp. $A^-(u)$)
is the set of arcs with $u$ as
its head (resp. tail), and 
$f$ is called a {\it nowhere-zero $\Gamma$-flow} of $D$
if it is a $\Gamma$-flow and $f(a)\ne 0$ 
for each $a\in A(D)$.

An example is shown in Figure~\ref{flow-f2}.

\begin{figure}[h!]
\centering

\begin{tikzpicture}[
	vertex/.style={circle,fill=black,draw=black,
		inner sep=0pt,minimum size=7.5pt,outer sep=0pt},
	every label/.style={font=\large,inner sep=2pt},
	directed edge/.style={
		-{Stealth[length=2.6mm,width=1.8mm]},
		line width=0.8pt,shorten <=0.5pt,shorten >=1pt
	},
	weight/.style={font=\large,inner sep=3pt}
	]
	\node[vertex,label=above:{$w$}] (w) at ( 0,    1.80) {};
	\node[vertex,label=left:{$x$}]  (x) at (-1.80, 0)    {};
	\node[vertex,label=below:{$y$}] (y) at ( 0,   -1.80) {};
	\node[vertex,label=right:{$z$}] (z) at ( 1.80, 0)    {};
	
	\draw[directed edge] (x) -- node[weight,midway,above left]  {$2$} (w);
	\draw[directed edge] (w) -- node[weight,midway,above right] {$1$} (z);
	\draw[directed edge] (w) -- node[weight,midway,right]       {$1$} (y);
	\draw[directed edge] (y) -- node[weight,midway,below left]  {$2$} (x);
	\draw[directed edge] (y) -- node[weight,midway,below right] {$3$} (z);
\end{tikzpicture}

\caption{A nowhere-zero $Z_4$-flow
\relabel{flow-f2}
}
\end{figure}

Nowhere-zero flows were introduced by Tutte \cite{tut2} 
as a dual concept to 
study proper colourings,
and the four-conjecture holds 
if and only if every bridgeless plane graph has a nowhere zero
$Z_4$-flow,
where $Z_p$ denotes 
the group of integers mod $p$.


An undirected graph $G$ 
is said to 
have a {\it nowhere-zero $\Gamma$-flow} if 
one of its orientations has such a flow. This definition is proper
as for any orientations $D_1$ and $D_2$,  it can be shown that 

\begin{center}
$D_1$ has a nowhere-zero $\Gamma$-flow
$\Longleftrightarrow $
$D_2$ has a nowhere-zero $\Gamma$-flow.
\end{center}

Furthermore,  
for any two Abelian groups $\Gamma_1$ and $\Gamma_2$ with the same order, 
$G$ has a nowhere-zero $\Gamma_1$-flow
if and only if it has 
a nowhere-zero $\Gamma_2$-flow.

For any positive integer $q$, 
a \red{\it nowhere-zero $q$-flow} is a nowhere-zero 
$Z$-flow $g$ such that $|g(a)| < q$ for all arcs $a$ in $D$.
A nowhere-zero 3-flow is shown  in Figure~\ref{flow-f3}.

\begin{figure}[htbp]
\centering

\begin{tikzpicture}[
	vertex/.style={circle,fill=black,draw=black,
		inner sep=0pt,minimum size=7.5pt,outer sep=0pt},
	every label/.style={font=\large,inner sep=2pt},
	directed edge/.style={
		-{Stealth[length=2.6mm,width=1.8mm]},
		line width=0.8pt,shorten <=0.5pt,shorten >=1pt
	},
	weight/.style={font=\large,inner sep=3pt}
	]
	\node[vertex,label=above:{$w$}] (w) at ( 0,    1.80) {};
	\node[vertex,label=left:{$x$}]  (x) at (-1.80, 0)    {};
	\node[vertex,label=below:{$y$}] (y) at ( 0,   -1.80) {};
	\node[vertex,label=right:{$z$}] (z) at ( 1.80, 0)    {};
	
	\draw[directed edge] (x) -- node[weight,midway,above left]  {$2$}  (w);
	\draw[directed edge] (w) -- node[weight,midway,above right] {$1$}  (z);
	\draw[directed edge] (w) -- node[weight,midway,right]       {$1$}  (y);
	\draw[directed edge] (y) -- node[weight,midway,below left]  {$2$}  (x);
	\draw[directed edge] (y) -- node[weight,midway,below right] {$-1$} (z);
\end{tikzpicture}

\caption{A nowhere-zero $3$-flow
\relabel{flow-f3}
}
\end{figure}

Tutte~\cite{tut2} showed 
that $G$ has a nowhere-zero $q$-flow
if and only if it has a nowhere-zero $Z_q$-flow.

 
 
 The following conclusion due to 
 Tutte shows that the $4$-color conjecture holds if and only if
 every plane graph has a 
 nowhere-zero $4$-flow.

\begin{theo}[Tutte \cite{tut1}]
	\relabel{4ct-flow}
A plane graph $G$ is $k$-face-colourable 
if and only if it has a nowhere-zero $k$-flow.
\end{theo}

\proof 
Let $G$ be a plane graph and $\stackrel{\rm \rightarrow }{G}$
be an orientation of $G$.

From a face colouring of $G$ with colouring $0,1,\cdots,k-1$,
we can get a nowhere-zero 
$k$-flow by assigning each arc the difference of 
the two values of its two sides: 
\red{the right-hand side to the arrow minus 
the other side.}
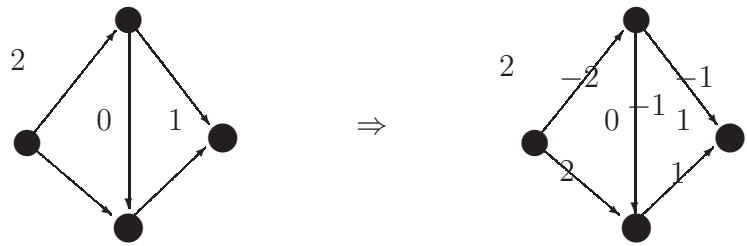
\proofend

An example of Theorem~\ref{4ct-flow} is shown 
in Figure~\ref{sec4-f3}. 

\begin{figure}[htp]
	\centering

	\begin{tikzpicture}[
		vertex/.style={circle,fill=black,draw=black,
			inner sep=0pt,minimum size=10pt,outer sep=0pt},
		directed edge/.style={
			-{Stealth[length=2.1mm,width=1.35mm]},
			line width=0.55pt,shorten <=0.5pt,shorten >=1pt
		},
		number/.style={font=\large,inner sep=2pt},
		edge label/.style={font=\large,fill=white,inner sep=1.5pt}
		]
		
		\def\DGPlaceVertices{%
			\node[vertex] (T) at ( 0,    1.70) {};
			\node[vertex] (L) at (-1.60, 0)    {};
			\node[vertex] (B) at ( 0,   -1.70) {};
			\node[vertex] (R) at ( 1.60, 0)    {};
		}
		
		\def\DGRegionLabels{%
			\node[number] at (-2.10,1.20) {$2$};
			\node[number] at (-0.50,0.08) {$0$};
			\node[number] at ( 0.56,0.08) {$1$};
		}
		
		\begin{scope}
			\DGPlaceVertices
			\draw[directed edge] (L) -- (T);
			\draw[directed edge] (T) -- (R);
			\draw[directed edge] (L) -- (B);
			\draw[directed edge] (B) -- (R);
			\draw[directed edge] (T) -- (B);
			\DGRegionLabels
		\end{scope}
		
		\node[font=\LARGE] at (3.60,0) {$\Rightarrow$};
		
		\begin{scope}[xshift=7.20cm]
			\DGPlaceVertices
			\draw[directed edge] (L)
			-- node[edge label,midway] {$-2$} (T);
			\draw[directed edge] (T)
			-- node[edge label,midway] {$-1$} (R);
			\draw[directed edge] (L)
			-- node[edge label,midway] {$2$} (B);
			\draw[directed edge] (B)
			-- node[edge label,midway] {$1$} (R);
			\draw[directed edge] (T)
			-- node[edge label,pos=0.31] {$-1$} (B);
			\DGRegionLabels
		\end{scope}
		
	\end{tikzpicture}

\caption{From a face coloring to a nonwhere-zero flow
\label{sec4-f3}
}

\end{figure}

Tutte's three flow conjectures
and Jaeger's weak 3-flow Conjecture are the main flow  problems studied by mathematicians, and we now present these problems below.

\begin{con}[Tutte's 5-flow Conjecture \cite{tut1}]
	
Every bridgless graph has a nowhere-zero $5$-flow.
\end{con}

\begin{con}[Tutte's 4-flow Conjecture \cite{tut66}]
	
Every bridgless graph with no Petersen minor 
has a nowhere-zero $4$-flow.
\end{con} 


\begin{con}[Tutte's 3-flow Conjecture (cf. \cite{Bondy})]
	
Every $4$-edge-connected graph has a 
nowhere-zero $3$-flow.
\end{con}


\begin{con}[Jaeger's weak 3-flow Conjecture \cite{jae1988}] 
There exists a fixed integer $k$ so that every $k$-edge-connected graph has a nowhere-zero $3$-flow.
\end{con} 

Thomassen~\cite{tho2012} proved  Jaeger's weak 3-flow conjecture for $k=8$, and 
his conclusion was extended to 
$k=6$ by Lov\'asz et al~\cite{lovasz2013}.

\subsection{Definition of flow polynomial
\label{sec4-2}
}

For any graph $G$ and any positive integer $t$, 
let $F(G, t)$ be \blue{the number of distinct nowhere-zero 
$Z_t$-flows of $G$}.
The function $F(G, t)$ is called the {\it flow polynomial} 
of $G$.

Tutte's 5-flow conjecture is equivalent to the statement that 
\red{$F(G, 5)>0$ for all bridgeless graph $G$}.

The flow polynomial $F(G,x)$ of a graph $G$
can be obtained 
recursively 
by the following  rules
(see Tutte~\cite{tut}):
\begin{equation}\relabel{eq1-2} 
F(G,x)=
\left \{
\begin{array}{ll}
1, &\mbox{if }E=\emptyset;\\
0, &\mbox{if }G\mbox{ has a bridge};\\
F(G_1,x)F(G_2,x),&\mbox{if }G=G_1\cup G_2;\\
(x-1)F(G\backslash e,x), &\mbox{if }e\mbox{ is a loop};\\
F(G / e,x)-F(G\backslash e,x),  \quad
&\mbox{if }e\mbox{ is not a loop nor a bridge},
\end{array}
\right.
\end{equation}
where $e\in E$.
The rules in (\ref{eq1-2}) can be 
verified by the definition of $F(G,x)$, 
and $G_1\cup G_2$ is the disjoint union of $G_1$ and $G_2$.

\begin{exa} 
\begin{enumerate}[(a)]
\item If $G$ is a cycle, then 
$F(G,x)=x-1$.
\item If $G=L_k$ is a graph with two distinct vertices
and $k$ edges joining them, then 
$$
F(G,x)=
\left ((x-1)^k+(-1)^k(x-1)\right )/x.
$$

\item 
$F(L_3,x)=(x-1)(x-2)$.

\item  If $G=K_4$, then 
$
F(G,x)=(x-1)(x-2)(x-3).
$
Thus $K_4$ has no nowhere-zero 3-flow.
\end{enumerate}
\end{exa} 

\subsection{Computation} 

By (\ref{eq1-2}), 
the following basic properties
for computing $F(G,x)$ can be 
established: 
\begin{enumerate}[(a)]
	\item If $G_1,G_2,\cdots, G_k$ are components of $G$,
	or if $G$ is connected with 
	blocks
	$G_1,G_2,\cdots, G_k$, 
	then 
	$$
	F(G,x)=\prod_{1\le i\le k}
	F(G_i,x).
	$$

	\item  If $G$ has a bridge, 
	then $F(G,x)=0$.
	
	\item If $N_G(w)=\{u,v\}$ for some $w\in V(G)$, then 
	$$
	F(G,x)=F((G-w)\cdot uv,x),
	$$
	where $(G-w)\cdot uv$ is the graph obtained from $G-w$ 
	by identifying $u$ and $v$.

	\item (Jackson 2007) 
	Let $G$ be a bridgeless connected graph, $v$ be a vertex of $G$, 
	$e = u_1u_2$ be an edge of $G$, and 
	$H_1$ and $H_2$ be edge-disjoint subgraphs of 
	$G$ such that $E(H_1) \cup E(H_2) = E(G \backslash e)$, 
	$V(H_1)\cap V(H_2) = \{v\}$, 
	$V(H_1)\cup V(H_2) =V(G)$,
	$u_1\in V (H_1)$ and $u_2\in V(H_2)$, as shown 
	as shown 
	in Figure~\ref{sec4-f4}.
	Then 
	$$
	F(G, x) =
	\frac{F(G_1, x)F(G_2, x)}{x -1}.
	$$
	where $G_i=H_i+vu_i$ 
	for $i\in \{1, 2\}$.
	
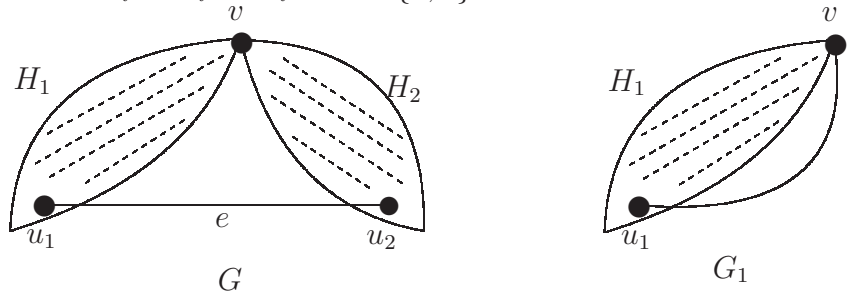
	\begin{figure}[h!]
		\centering

\begin{tikzpicture}[
	line cap=round,
	line join=round,
	boundary/.style={line width=0.35pt},
	graph edge/.style={line width=0.50pt},
	hatching/.style={line width=0.25pt,dash pattern=on 2.5pt off 2.2pt},
	vertex/.style={circle,draw=black,fill=black,
		inner sep=0pt,minimum size=7.5pt},
	vertex label/.style={font=\large,inner sep=2pt},
	region label/.style={font=\Large,inner sep=2pt},
	graph label/.style={font=\Large,inner sep=2pt}
	]
	
	\def\HGCoordinates{%
		\coordinate (v)     at ( 0,    2.60);
		\coordinate (ltip)  at (-3.20, 0);
		\coordinate (rtip)  at ( 2.55, 0);
		\coordinate (u1)    at (-2.72, 0.36);
		\coordinate (u2)    at ( 2.05, 0.36);
	}
	
	\def\HGLeftPath{%
		(v)
		.. controls (-2.05,2.73) and (-3.22,1.81) .. (ltip)
		.. controls (-1.25,0.57) and (-0.30,1.39) .. (v)
		-- cycle
	}
	\def\HGRightPath{%
		(v)
		.. controls (1.93,2.72) and (2.61,1.70) .. (rtip)
		.. controls (1.18,0.25) and (0.28,1.24) .. (v)
		-- cycle
	}
	
	\def\HGDrawLeft{%
		\begin{scope}
			\clip \HGLeftPath;
			\foreach \k in {0,...,12} {
				\draw[hatching]
				(-3.5,{-1.96+0.32*\k}) -- (0.3,{0.168+0.32*\k});
			}
		\end{scope}
		\draw[boundary] \HGLeftPath;
	}
	\def\HGDrawRight{%
		\begin{scope}
			\clip \HGRightPath;
			\foreach \k in {0,...,12} {
				\draw[hatching]
				(-0.3,{0.216+0.32*\k}) -- (2.9,{-2.088+0.32*\k});
			}
		\end{scope}
		\draw[boundary] \HGRightPath;
	}
	
	\begin{scope}
		\HGCoordinates
		\HGDrawLeft
		\HGDrawRight
		
		\draw[graph edge] (u1)
		-- node[font=\large,midway,below,inner sep=3pt] {$e$} (u2);
		
		\foreach \p in {v,u1,u2}
		\node[vertex] at (\p) {};
		
		\node[vertex label,anchor=south] at (0,2.80) {$v$};
		\node[vertex label,anchor=north] at (-2.72,0.12) {$u_1$};
		\node[vertex label,anchor=north] at (2.05,0.12) {$u_2$};
		\node[region label] at (-2.96,2.12) {$H_1$};
		\node[region label] at (2.32,2.08) {$H_2$};
		\node[graph label] at (-0.16,-0.65) {$G$};
	\end{scope}
	
	\begin{scope}[xshift=8.10cm]
		\HGCoordinates
		\HGDrawLeft
		
		\draw[graph edge] (u1)
		.. controls (-0.52,0.19) and (0.22,0.80) .. (v);
		
		\foreach \p in {v,u1}
		\node[vertex] at (\p) {};
		
		\node[vertex label,anchor=south] at (0,2.80) {$v$};
		\node[vertex label,anchor=north] at (-2.72,0.12) {$u_1$};
		\node[region label] at (-2.96,2.12) {$H_1$};
		\node[graph label] at (-1.40,-0.65) {$G_1$};
	\end{scope}
	
\end{tikzpicture}

		\caption{$G \backslash e$ 
			has a cut-vertex $v$
		\label{sec4-f4}
	}
		
	\end{figure} 
	
	\item  (Jackson  \cite{jac2007})  
	Let $G$ be a bridgeless connected graph, 
	$S$ be a $2$-edge-cut of $G$, 
	and $H_1$ and $H_2$ be the sides of $S$, as shown 
	in Figure~\ref{sec4-f5}.
	Let $G_i$ be obtained from $G$ by contracting 
	$E(H_{3-i})$, for $i\in \{1, 2\}$. Then
	$$
	F(G, x) =
	\frac{F(G_1, x)F(G_2, x)}{x -1}.
	$$
	
	
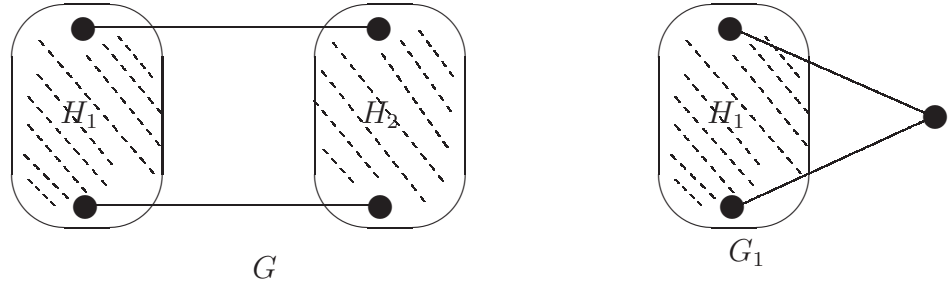
\begin{figure}[h!]
		\centering

		\begin{tikzpicture}[
			line cap=round,
			line join=round,
			boundary/.style={line width=0.4pt,rounded corners=18pt},
			hatching/.style={line width=0.25pt,dash pattern=on 2.5pt off 2.2pt},
			graph edge/.style={line width=0.5pt},
			vertex/.style={circle,fill=black,draw=black,
				inner sep=0pt,minimum size=8pt},
			region label/.style={font=\Large,fill=white,inner sep=3pt},
			graph label/.style={font=\Large,inner sep=2pt}
			]
			
			\def\DrawRoundedRegion#1{%
				\begin{scope}
					\clip[rounded corners=15pt]
					(-0.82,-1.27) rectangle (0.82,1.27);
					\foreach \k in {-10,...,10} {
						\draw[hatching]
						(-1.2,{1.2+0.34*\k}) -- (1.2,{-1.2+0.34*\k});
					}
				\end{scope}
				\draw[boundary] (-0.95,-1.40) rectangle (0.95,1.40);
				\node[region label] at (0,0) {$#1$};
			}
			
			\begin{scope}
				\DrawRoundedRegion{H_1}
				\begin{scope}[xshift=3.70cm]
					\DrawRoundedRegion{H_2}
				\end{scope}
				
				\coordinate (a) at (0,1.10);
				\coordinate (b) at (0,-1.10);
				\coordinate (c) at (3.70,1.10);
				\coordinate (d) at (3.70,-1.10);
				
				\draw[graph edge] (a) -- (c);
				\draw[graph edge] (b) -- (d);
				
				\foreach \p in {a,b,c,d}
				\node[vertex] at (\p) {};
				
				\node[graph label] at (1.85,-1.95) {$G$};
			\end{scope}
			
			\begin{scope}[xshift=8.00cm]
				\DrawRoundedRegion{H_1}
				
				\coordinate (a) at (0,1.10);
				\coordinate (b) at (0,-1.10);
				\coordinate (c) at (2.55,0);
				
				\draw[graph edge] (a) -- (c) -- (b);
				
				\foreach \p in {a,b,c}
				\node[vertex] at (\p) {};
				
				\node[graph label] at (0.80,-1.95) {$G_1$};
			\end{scope}
			
		\end{tikzpicture}

		\caption{$G$ has a $2$-edge-cut
		\label{sec4-f5}
	}
	\end{figure} 
	
	\item (Jackson \cite{jac2007})  
	Let $G$ be a bridgeless connected  graph, 
	$S$ be a $3$-edge-cut of $G$, 
	and $H_1$ and $H_2$ be the sides of $S$. 
	Let $G_i$ be obtained from $G$ by contracting 
	$E(H_{3-i})$, for $i\in \{1, 2\}$.
	Then
	$$
	F(G, x) = 
	\frac{F(G_1, x)F(G_2, x)}
	{(x -1)(x -2)}. 
	$$
\end{enumerate}

\subsection{Expansion}

For a connected graph $G=(V,E)$,   
the following expression shows that $F(G,x)$ is indeed a polynomial in $x$:
\begin{equation}
	\relabel{flow-eq4}
	F(G,x)=\sum\limits_{A\subseteq E}
	(-1)^{|E|-|A|}x^{|A|-|V|+c(A)},
\end{equation}
where 
$c(A)$ is the number of components 
in the subgraph $(V,A)$.
(\ref{flow-eq4}) can also be 
proved by applying (\ref{eq1-2}).

Assume that $G$ is 
a bridgeless and connected
graph of order $n$ and size $m$.
By (\ref{flow-eq4}), 
$F(G,x)$ can be expressed as
follows: 
\equ{sec4-e4}
{
F(G,x)=x^{m-n+1}-b_1x^{m-n}+b_2x^{m-n-1}+\cdots+(-1)^{m-n+1}b_{m-n+1},
}
where 
\begin{enumerate}[(a)]
	\item for any $1\le i\le m-n+1$, 
	$b_i$ is positive
	(can be proved by induction) and 
	$$
	b_{i}=\nu_{i,1}-\nu_{i+1,2}+\nu_{i+2,3}-\dots
	=\sum_{j\ge 1}(-1)^{j-1}\nu_{i+j-1,j},
	$$
	where 
$\nu_{k,j}$ is the number of subsets $A$ of $E$ 
such that $|A|=k$ and $c(G-A)=j$; and 

\item 
If $G$ has no $2$-edge-cut,
then $b_{1}=m$ and $b_{2}={m\choose 2}-\gamma$,
where $\gamma$ is the number of $3$-edge-cuts of $G$.
\end{enumerate}

\subsection{Connection to other polynomials}

By (\ref{sec2-e1}) and 
(\ref{flow-eq4}), the following relation between $F(G,x)$ and 
the Tutte polynomial $T_G(x,y)$ of $G=(V,E)$ can be 
established:
\begin{equation}\relabel{flow-eq3}
F(G,x)=(-1)^{|E|-|V|+c(E)}
T_G(0, 1-x),
\end{equation}
where $c(E')$ is the number of components 
in the spanning subgraph $(V,E')$ of $G$.
The chromatic polynomial 
$\chi(G,x)$
is considered to be a dual  polynomial of 
$F(G,x)$, 
as for any plane graph $G$, 
\equ{sec4-e5}
{
	\chi(G,x)=x
	F(G^{*},x), 
}
where $G^*$ is  the dual of a plane  graph $G$.
This conclusion can be obtained 
by applying (\ref{flow-eq3}), 
the fact $T_G(x,y)=T_{G^*}(y,x)$ (see 
Proposition~\ref{tutte-pro2})
and the relation between 
$T_G(x,y)$ and $\chi(G,x)$
(see page \pageref{spe-case}).

\subsection{Known results}

In this subsection, we present some 
known results on the
locations of roots of $F(G,x)$.

\begin{theo}
[Waklin \cite{wak}] 
\label{sec4-th1}
Let $G=(V,E)$ be a bridgeless connected graph
with block number $b$. Then
\begin{enumerate}[(a)]
\item
$F(G, x)$ is non-zero with sign 
$(-1)^{|E|-|V|+1}$ for $x\in (-\infty, 1)$;
\item
$F(G, x)$ has a zero of multiplicity 
$b$ at $x = 1$; and 
\item  
$F(G, x)$ is non-zero with sign 
$(-1)^{|E|-|V|+b-1}$ for $x\in (1, 32/27]$.
\end{enumerate}
\end{theo}

\begin{theo}
	[Jackson \cite{jac2007}] 
	\label{sec4-th6}
	Let $G$  be a 3-connected cubic graph with $n$ vertices and
	$m$ edges. Then
	\begin{enumerate}[(a)]
		\item $F(G, x)$ is non-zero 
		with sign $(-1)^{m-n}$
		in the interval $(1,2)$;
		\item  $F(G, x)$  
		has a zero of multiplicity $1$ at $x = 2$;
		\item $F(G, x)$ is non-zero with sign 
		$(-1)^{m-n+1}$ for $x\in (2, d)$, 
		where $d\approx 2.546$ is the flow root of the cube in 
		$(2, 3)$.
	\end{enumerate}
\end{theo}

Jackson~\cite{jac2} showed that for any bridgeless graph $G$
of order $n$, 
$F(G,x)>0$ holds for all real $x$ 
with
$x\ge 2\log_2 n$.

\begin{theo}[Jackson \cite{jac2007}] 
	\label{sec4-th2}
If $G$ has at most one 
vertex of degree larger than $3$, then 
$F(G, x)$ is non-zero in the interval $(1,2)$.
\end{theo}

\begin{theo}
[Dong  \cite{don2015b, don2015}]
\label{sec4-th3}
If $G$ has at most two
vertices of degrees larger than $3$, then 
$F(G, x)$ is non-zero in the interval $(1,2)$.
 More generally, if all vertices in 
the set $W:=\{u\in V(G): d(u)\ge 4\}$
are dominated by one component of $G-W$, 
then $F(G, x)$ is non-zero in the interval $(1,2)$.
\end{theo}

\begin{theo}
[Kung and Royle \cite{kung2011}] 
\label{sec4-th4}
If $G$ is a bridgeless graph, 
then its flow roots are integral if and only
if $G$ is the dual of a chordal and plane graph.
\end{theo}

\begin{theo}
[Dong \cite{don2018}] 
\label{sec4-th5}
For any  bridgeless graph $G$,
if $F(G,x)$ has real roots only,
then either every root of $F(G,x)$ is in the set $\{1,2,3\}$ 
or $F(G,x)$ has at least 9 roots in $(1,2)$.
\end{theo}

By Theorem~\ref{sec4-th1} (b), 
Theorem~\ref{sec4-th6} (b), 
and Theorems~\ref{sec4-th3} 
and~\ref{sec4-th5}, the following conclusion holds.

\begin{cor}
	\label{sec4-cor1}
For any $2$-connected graph $G$, if $F(G,x)$ is non-zero in $(1,2)$ and 
$F(G,x)$ has real roots only,  then 
\equ{sec4-e6}
{
F(G,x)=(x-1)(x-2)^{n_2}
(x-3)^{n_3}
}
for some non-negative 
integers $ n_2$ and $n_3$,
where $n_2=1$ when $G$ is 
$3$-connected.
\end{cor} 

More details on the roots 
of $F(G,x)$ 
can be found in a survey due to
Dong \cite{don2019}.


\subsection{Open problems}

The following problem due to Welsh was proposed in 1970s.

\begin{con}[Welsh \cite{wel}]
	\relabel{welsh-con4}
For any bridgeless graph $G$, $F(G,x)>0$ for all real numbers 
$x\in (4,\infty)$.
\end{con}

\vspace{0.5 cm}

Haggard, Pearce and Royle \cite{hag}
showed that the generalised Petersen graph $G_{16,6}$  
has real flow roots around $4.0252205$ and $4.2331455$, 
where the generalized Petersen graph $G_{n,k}$ 
for $n\ge 3$ and $1\le k\le \lfloor (n-1)/2\rfloor$ 
is the graph with vertex set $\{u_i,v_i: 1\le i\le n\}$ 
and edge set 
$\{u_iv_i, u_iu_{i+1}, v_iv_{i+k}: 1\le i\le n\}$,
where $v_s$ for $s>n$ is considered as $v_t$,
where $t$ is the integer with $1\le t\le n$ 
such that $s-t$ is a multiple of $n$. They then revised Welsh's conjecture as follows.

\begin{con}[Haggard, Pearce and Royle \cite{hag}]
\relabel{hag-con5}
For any bridgeless graph $G$, $F(G,x)>0$ for all real numbers 
$x\in [5,\infty)$.
\end{con}

Conjecture~\ref{hag-con5} was recently disproved by 
Jacobsen and Salas~\cite{jaco} who found counter-examples 
by studying the subfamily of 
generalised Petersen graphs $G_{nr,r}$ for $n\ge 2$ and $r\ge 2$. 

\begin{theo}[Jacobsen and Salas \cite{jaco}]
	\relabel{jaco-them1}
The value $q=5$ is an isolated accumulation point of real zeros of the flow polynomial $F(G,q)$ for the families of bridgeless graphs $G_{6n,6}$ and $G_{7n,7}$ with $n\ge 3$. Moreover:
\begin{enumerate}[(a)]
\item 
There is a sequence of real zeros $\{q_n\}$ of the flow polynomials $F(G_{6n,6},q)$ that converges to $q=5$ from below.

\item
There is a sequence of real zeros $\{q_n\}$ of the flow polynomials $F(G_{7n,6},q)$
that converges to $q=5$.
The sub-sequence with odd (resp.even) $n$ converges to $q=5$ from above(resp. below).

\item
The bridgeless graph $G_{119,7}$
has flow roots at $q\approx  5.00002$ and 
$q\approx  5.16534$ (where $\approx $ means ``within $10^{-5}$").
\item
The value $q'\approx 5.235261$ (where $\approx $ means ``within $10^{-6}$")
is an accumulation point of real zeros of the flow polynomials 
$F(G_{7n,7},q)$.
In particular,
the sub-sequence for odd $n$ of the real zeros $\{q_n\}$ of the flow polynomials $F(G_{7n,7},q)$ converges to $q'$
from below.
\end{enumerate}
\end{theo}

\begin{con}[Jacobsen and Salas \cite{jaco}]\relabel{jaco-con6}
For any bridgeless graph $G$, $F(G,q)>0$ for all real numbers 
$q\in [6,\infty)$.
\end{con}

\begin{con}
[Dong \cite{don2018,don2020}] \relabel{dong-con2}
For any bridgeless graph $G$, 
if $F(G,q)$ has real roots only, 
then all roots of $F(G,q)$ are integral.
\end{con}

\section{$\sigma$-polynomial $\sigma(G,x)$}
\relabel{sect-sigma}

\subsection{Definition}

For a graph $G$ of order $p$,  its chromatic 
polynomial can be written as 
\equ{sec6-e1}
{
\chi(G,x)=\sum_{0\le i\le p} a_i(G)\cdot (x)_i,
}
where $(x)_i=x(x-1)\cdots (x-i+1)$.
Obviously, 
$a_i(G)=0$ for all $i<\chi(G)$, and 
$a_i(G)$ is positive integer for $\chi(G)\le i\le p$.
Actually $a_i(G)$ is the number of partitions 
of $V(G)$ into $i$ non-empty independent sets.

For any non-adjacent pair of vertices $u$ and $v$ in $G$,
by (\ref{sec5-e1}), 
\equ{sec6-e2}
{
a_i(G)=a_i(G+uv)+a_i(G\cdot uv),
}
where $G+uv$ is the graph obtained from $G$ by adding a new 
edge joining $u$ and $v$ and 
$G\cdot uv$ is the graph obtained from $G$ 
by identifying $u$ and $v$,
and removing all but one parallel edges between any two vertices.

Now we define two functions 
${\sigma}(G,x)$
and $\bar {\sigma}(G,x)$
as follows:
\equ{sec6-e3}
{
\sigma(G,x)=\sum_{0\le i\le p} a_i(G) x^i,
\quad 
\bar {\sigma}(G,x)=\sum_{0\le i\le p}  i!a_i(G) x^i.
}

$\sigma(G,x)$ was defined by Korfhage~\cite{kor1978}, 
although the original function he introduced was 
actually $\sigma(G,x)/x^{\chi(G)}$.

The {\it adjoint polynomial} $h(G,x)$ defined below
was introduced by Liu Ruying~\cite{liu1987}:
\equ{sec6-e4}
{
h(G,x)=\sum_{i} h_i x^i,
}
where $h_i$ is the number of partitions of $V(G)$ 
into exactly $i$ non-empty subsets each of which is a clique of $G$. 
Observe that $h(G,x)=\sigma(\bar G,x)$, where $\bar G$ is the complement of $G$.

\subsection{Examples}

Let $K_n, N_n$, $P_n$ and $C_n$ denote 
the complete graph,
the null graph, 
the path graph 
and the cycle graph of order $n$, respectively. 

By the definition of $\sigma(G,x)$ given in (\ref{sec6-e3}), it is not difficult to obtain the following conclusions on $\sigma(G,x)$ for some special graphs $G$.

\begin{enumerate}[(1)]
\item $\sigma(K_n,x)=x^n,$ as $\chi(K_n,x)=(x)_n$.

\item $\sigma(N_n,x)=
\sum\limits_{1\le k\le n}S(n,k)x^k$,
as
\equ{sec6-e6}
{
\chi(N_n,x)=x^n
=\sum_{1\le k\le n}
S(n,k)(x)_k.
}

\item If $G$ is the complete $r$-partite graph 
$K_{m_1,m_2,\cdots,m_r}$,
then 
\equ{sec6-e7}
{
\sigma(G,x)=
\prod_{i=1}^r 
\left ( 
\sum_{1\le k\le m_i}
S(m_i,k)x^k
\right ).
}
This result can be deduced 
from (\ref{chap5-factor5-eq1}).

\item  (Liu~\cite{liu1987})
$
\sigma(\bar P_n, x)=
\sum\limits_{i\le n}{i\choose n-i}x^i.
$

\item (Dong et al \cite{don2002})
\equ{sec6-e10}
{
\sigma(\bar P_n, x)=x^{\lceil n/2\rceil }
\prod_{s=1}^{\lfloor n/2\rfloor}
\left (x+2+2\cos \frac{2s\pi}{n+1}\right )
}
and 
\equ{sec6-e11}
{
\sigma(\bar C_n, x)=x^{\lceil n/2\rceil }
\prod_{s=1}^{\lfloor n/2\rfloor}
\left (x+2+2\cos \frac{(2s-1)\pi}{n+1}\right ).
}

\item (Brenti~\cite{bre1992})
$\sigma(G)=
x^n(x+1)^n$
for $G\cong K_{\underbrace{2,2,\cdots,2}_n}$.
\end{enumerate}

\subsection{Computation}

Recall that the {\it joint} of disjoint graphs $G_1$ and $G_2$, 
denoted by $G_1+ G_2$, 
is obtained from the disjoint union of 
graphs $G_1$ and $G_2$ by adding edges 
joining each $u\in V(G_1)$ to each $v\in V(G_2)$.
Thus, for any $k\ge 1$, 
\equ{sec6-e12}
{
a_k(G_1+ G_2)=\sum_{i+j=k}a_i(G_1)a_j(G_2)
}
and
\equ{sec6-e13}
{
\sigma(G_1+ G_2,x)
=\sigma(G_1,x)\sigma(G_2,x).
}
If $u$ and $v$ are non-adjacent vertices in $G$, 
by (\ref{sec6-e2}), 
\equ{sec6-e14}
{
\sigma(G,x)=\sigma(G+uv,x)+\sigma(G\cdot uv,x).
}
Hence the function $\sigma(G,x)$ can be computed by the following recursive
formula:
\begin{equation}
	\label{sec6-e15}
	\sigma(G,x)=
	\left \{
	\begin{array}{ll}
		x^{n}, &\mbox{if }
		G\cong K_n;\\
		\sigma(G_1,x)\sigma(G_2,x),
		&\mbox{if }G=G_1+ G_2;\\
		\sigma(G+uv,x)
		+\sigma(G\cdot uv, x),  \quad
		&\mbox{if }u,v\in V, 
		uv\not\in E(G).
	\end{array}
	\right.
\end{equation}

For a graph $G$ of order $p$ and size $q$, 
some leading coefficients of $\sigma(G,x)$ can be determined directly, 
as shown below:
\begin{enumerate}[(a)]
\item $a_p(G)=1$;

\item $a_{p-1}={p\choose 2}-q$;
and 

\item  (Brenti \cite{bre1992})
$
a_{p-2}={q\choose 2}-q{q-1\choose 2}
+{p\choose 3}{3p-5\choose 4}-t(G),
$
where $t(G)$ is the number of triangles in $G$.
\end{enumerate} 

\subsection{Open problems}

A graph $G$ is said to be {\it $\sigma$-real}
if $\sigma(G,x)$ has real zeros only, and 
it is said to be {\it $\sigma$-unreal}
if it is not $\sigma$-real.

The $w$-unreal and $\tau$-unreal graphs 
are defined similarly with respect to 
$w$-polynomial and $\tau$-polynomial respectively.
These polynomials will be introduced
Section~\ref{sect-w} and Section~\ref{sect-tau} respectively.

\begin{theo}[Brenti \cite{bre1992}]
	\relabel{sec6-th1}
$G$ is $\sigma$-real  
if one of the following conditions is satisfied:
\begin{enumerate}[(a)]
\item $\bar G$ is a comparability graph, where a graph $H$ is called 
a {\it comparability graph} if there exists a partial order $\preceq $
such that $uv\in E(H)$ if and only if $u\ne v$ and $u\preceq v$
or $v\preceq u$;

\item $\chi(G)\ge |V(G)|-2$; 

\item $\bar G$ is $K_3$-free;

\item there exists a simplicial vertex $u$ in $G$ 
such that $G-u$ is $\sigma$-real;

\item $G=G_1+G_2$, where 
both $G_1$ and $G_2$
are $\sigma$-real;

\item $G=G_1\cup G_2$\footnote{
	Here $G_1$ and $G_2$ 
	are not vertex-disjoint,
	and for any $u,v\in V(G_1)\cap V(G_2)$, 
	$uv\in E(G_1)$ if and only if $uv\in E(G_2)$.
$G_1\cup G_2$ is 
the graph with vertex set $V(G_1)\cup V(G_2)$ 
and edge set $E(G_1)\cup E(G_2)$.}, 
where each $G_i$ is $\sigma$-real
and $G_1\cap G_2$ is complete; 

\item $\bar\sigma(G,x)$ has real zeros only; and 

\item $w(G,x)$ has real zeros only.
\end{enumerate}
\end{theo}

 The $\sigma$-unreal, $w$-unreal and $\tau$-unreal 
connected graphs on up to 9 vertices
were determined by 
Cameron, Colbourn, Read and Wormald~\cite{cam},
and the numbers of $\sigma$-unreal,  $w$-unreal 
and $\tau$-unreal connected graphs of orders 
from $3$ to $9$ are shown below:
\begin{center}
\begin{tabular}{c|ccccccc}
order & 3 & 4& 5&6&7&8&9 \\ \hline
no. $\sigma$-unreal con. graphs 
& 0 & 0 &0 &0 &0 &2 & 42 \\ \hline
no. $w$-unreal con. graphs 
& 0 & 1 & 3 & 16 & 116 & 1237 & 22515  \\ \hline
no. $\tau$-unreal con. graphs 
& 0 & 0 &0 &0 &0 &0 & 0 \\ 
\end{tabular}
\end{center}
Note that, in the above table, 
$w$-unreal (resp. $\tau$-unreal) con. graphs 
refer to connected graphs whose $w$-polynomials 
(resp. $\tau$-polynomials) 
have unreal roots.  

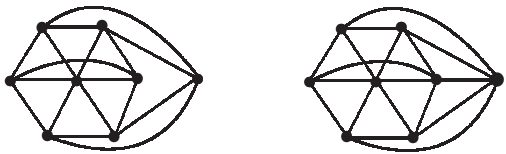
\begin{figure}[h!]
\begin{center}
\begin{tikzpicture}[
	line width=0.55pt,
	line cap=round,
	line join=round,
	vertex/.style={circle,draw=black,fill=black,
		inner sep=0pt,minimum size=6pt},
	graph label/.style={font=\Large,inner sep=2pt}
	]
	
	\def\GGCoordinates{%
		\coordinate (L) at (-1.24, 0);
		\coordinate (A) at (-0.62, 1.05);
		\coordinate (B) at ( 0.62, 1.05);
		\coordinate (C) at ( 0,    0);
		\coordinate (D) at ( 1.24, 0);
		\coordinate (E) at (-0.62,-1.05);
		\coordinate (F) at ( 0.62,-1.05);
		\coordinate (R) at ( 2.45, 0);
	}
	
	\def\GGCommonEdges{%
		\draw (L) -- (A) -- (B) -- (D) -- (F) -- (E) -- cycle;
		
		\foreach \p in {L,A,B,D,F,E}
		\draw (C) -- (\p);
		
		\draw (B) -- (R) -- (F);
		
		\draw (A) .. controls (0.30,2.00) and (1.55,1.10) .. (R);
		\draw (E) .. controls (0.55,-1.50) and (1.75,-1.52) .. (R);
		\draw (L) .. controls (-0.58,0.52) and (0.58,0.52) .. (D);
	}
	
	\def\GGDrawVertices{%
		\foreach \p in {L,A,B,C,D,E,F,R}
		\node[vertex] at (\p) {};
	}
	
	\begin{scope}
		\GGCoordinates
		\GGCommonEdges
		\GGDrawVertices
		\node[graph label] at (0.60,-1.92) {$G_1$};
	\end{scope}
	
	\begin{scope}[xshift=5.85cm]
		\GGCoordinates
		\GGCommonEdges
		\draw (D) -- (R);
		\GGDrawVertices
		\node[graph label] at (0.60,-1.92) {$G_2$};
	\end{scope}
	
\end{tikzpicture}

\end{center}

\caption{$\sigma(G_1)=x^8 + 11x^7 + 38x^6 + 36x^5 + 11x^4 +x^3$
and $\sigma(G_2)=x^8 + 10x^7 + 30x^6 + 31x^5 + 10x^4 +x^3$,
and both polynomials have 
non-real zeros \cite{bre1994}.
\label{sec5-f1} 
}

\end{figure}

The two $\sigma$-unreal connected graphs on $8$
vertices are shown in Figure~\ref{sec5-f1}
(see \cite{bre1994}),
and 
the $\sigma$-unreal connected graphs on $9$ 
vertices are listed in~\cite{bre1994}.

\begin{con}[Brenti \cite{bre1992}] \relabel{con2-1}
Let $G$ be a simple graph of order $p$.
Then $G$ is $\sigma$-real if $\chi(G)\ge p-3$.
\end{con}

\section{$w$-polynomial $w(G,x)$}\relabel{sect-w}

\subsection{Definition} 
Let $G$ be a graph of order $p$. 
For any $1\le i\le p$, let $w_i(G)$
(or simply $w_i$) be the number 
determined by the following 
expression:
 \equ{sec7-e1}
 {
 	\chi(G,x)=\sum_{0\le i\le p} w_i{x+p-i\choose p}.
 }
By Theorem~\ref{theo1-2}, 
$w_i$ is the number of ordered pairs $(D,\pi)$,
where $D\in \seta(G)$, $\seta(G)$ is the set of 
acyclic orientations of $G$
and $\pi$ is an 
order-preserved permutation in $\OP(D)$ 
with $\rho(\pi)=i-1$.

The polynomial $w(G,x)$ is defined as follows (see \cite{bre1992}):
\equ{sec7-e2}
{
w(G,x)=\sum_{0\le i\le p} w_ix^i.
}

\begin{exa}\label{sec7-exa1}
As
$
\chi(K_p,x)
=(x)_p=p! {x\choose p},
$
we have 
$
w(K_p,x)=p!x^p.
$
\end{exa}

\subsection{Computation}

By Example~\ref{sec7-exa1}, 
$w(G,x)$ can be determined if $G$ is complete.
If $G$ is not complete, we will 
apply Lemma~\ref{lem3-2} below.

\begin{lem}\relabel{lem3-1}
If $u$ and $v$ are 
non-adjacent vertices in 
$G$, then 
\equ{sec7-e3}
{
w_i(G)=w_i(G+uv)+w_i(G\cdot uv)-w_{i-1}(G\cdot uv).
}
\end{lem}

\proof By (\ref{sec5-e1}), 
\equ{sec7-e4}
{
\chi(G,x)=\chi(G+uv,x)+\chi(G\cdot uv,x)
}
and 
\begin{eqnarray*}
&&\chi(G\cdot uv,x)\\
&=&\sum_{0\le i\le p-1}
w_i(G\cdot uv) {x+p-1-i\choose p-1}\\
&=& \sum_{0\le i\le p-1}
w_i(G\cdot uv){x+p-i\choose p}
-\sum_{0\le i\le p-1}
w_i(G\cdot uv){x+p-1-i\choose p}\\
&=& \sum_{0\le i\le p-1}
w_i(G\cdot uv){x+p-i\choose p}
-\sum_{1\le j\le p}
w_{j-1}(G\cdot uv){x+p-j\choose p}.
\end{eqnarray*}
Thus (\ref{sec7-e3}) follows.
\proofend

By the definition of $w(G,x)$ and  Lemma~\ref{lem3-1}, 
the next result follows directly. 

\begin{lem}\relabel{lem3-2}
If $u,v$ are not adjacent in $G$, then 
$$
w(G,x)=w(G+uv,x)+(1-x)w(G\cdot uv,x).
$$
\end{lem}

By Example~\ref{sec7-exa1}
and Lemma~\ref{lem3-2}, 
the polynomial $w(G,x)$ for 
any graph $G=(V,E)$ can be determined by the following 
recursive formula:
\begin{equation}
	\label{sec7-e8}
	w(G,x)=
	\left \{
	\begin{array}{ll}
		p!x^p, &\mbox{if }
		G\cong K_p;\\
	w(G+uv,x)+(1-x)w(G\cdot uv,x),
		\quad &\mbox{if } u,v\in V, 
		uv\not\in E.
	\end{array}
	\right.
\end{equation}

\begin{exa}
Applying (\ref{sec7-e8}) yields that 
$$
w(P_3,x)=w(K_3,x)+(1-x)w(K_2,x)=3!x^3+(1-x)2!x^2
=4x^3+2x^2
$$
and 
\begin{eqnarray*}
w(N_3,x)
&=&w(K_2\cup K_1,x)+(1-x)w(N_2,x)\\
&=&w(P_3,x)+(1-x)w(K_2,x)+(1-x)w(K_2,x)+(1-x)^2w(K_1,x)\\
&=&4x^3+2x^2+2(1-x)2!x^2+(1-x)^2x\\
&=&x^3+4x^2+x,
\end{eqnarray*}
where $K_2\cup K_1$ is the disjoint union of $K_2$ and $K_1$.
\end{exa}

\subsection{Properties}

\begin{theo}[Brenti\cite{bre1992}] 
\relabel{th3-1} 
For any graph $G$ of order $p$, 
\equ{sec7-e5}
{
\sum_{i\ge 0} \chi(G,i) x^i = 
\frac{w(G,x)}{(1-x)^{p+1}}.
}
\end{theo}

\proof Observe that 
\equ{sec7-e6}
{
\frac{w(G,x)}{(1-x)^{p+1}}
=(w_0+w_1x+\cdots+w_px^p)\sum_{j\ge 0}{p+j\choose p}x^j.
}
Thus, the coefficient of $x^i$ 
in the right-hand side of 
(\ref{sec7-e6}) 
is 
$$
\sum_{0\le j\le i}w_{i-j}{p+j\choose p}
=\sum_{0\le k\le i}w_{k}{p+i-k\choose p}
=\sum_{0\le k\le p}w_{k}{p+i-k\choose p}
=\chi(G,i).
$$
Thus the result follows.
\proofend

\begin{pro}[Brenti \cite{bre1992}]\relabel{pro3-1}
Let $G$ be a graph of order $p$. Then 
\begin{enumerate}[(a)]
\item $w_i=0$ for $i<\chi(G)$;

\item $w(G,1)=\sum_{i\le p} w_i=p!$;

\item $w_i$ is positive for $\chi(G)\le i\le p$; and 

\item $w_p$ is the number of acyclic orientations of $G$.
\end{enumerate} 
\end{pro}

\proof (a) By Theorem~\ref{th3-1}, 
$$
w(G,x)=(1-x)^p\sum_{i\ge \chi(G)}\chi(G,i)x^i,
$$
implying that (a) holds.
It also follows from Lemma~\ref{lem3-1}.

(b) It holds when $G$ is $K_p$, as 
$$
w(K_p,x)=p!x^p.
$$
Then, by Lemma~\ref{lem3-2}, $w(G,1)=p!$ for any graph $G$ of order $p$.

(c) It directly follows from Theorem~\ref{theo1-2}.

(d) Taking $x=-1$  yields that 
$$
\chi(G,-1)=\sum_{0\le i\le p} w_i{-1+p-i\choose p}
=w_p{-1\choose p}=(-1)^p w_p.
$$
As $(-1)^p\chi(G,-1)$ is the number of 
acyclic orientations of $G$, (d) holds.
\proofend

\begin{theo}[Brenti \cite{bre1992}] \relabel{w-sigmabar}
For any graph $G$ of order $p$, 
\equ{sec7-e7}
{
w(G,x)=(1-x)^p \bar \sigma \left (G, \frac{x}{1-x}\right ).
}
\end{theo}

\proof Let $z=x/(1-x)$, i.e., $x=z/(1+z)$. Then
the  identity 
(\ref{sec7-e7})
is equivalent to the following:
$$
w(G,z/(1+z))
=(1+z)^{-p} \bar \sigma \left (G, z\right )
$$
i.e., 
$$
\sum_{i\le p}w_iz^i(1+z)^{p-i}
=\sum_{0\le i\le p} i! a_i z^i,
$$
implying that 
\begin{equation}\relabel{eq3-8}
k!a_k=\sum_{i\le k}w_i{p-i\choose k-i}=\sum_{i\le k}w_i{p-i\choose p-k},
\quad \forall k\le p.
\end{equation}
By definition,
$$
\chi(G,x)=\sum_{i}w_i {x+p-i\choose p}
=\sum_{k}a_k(x)_k.
$$
As 
$$
{x+p-i\choose p}=\sum_{i\le k\le p}{x\choose k}{p-i\choose p-k}
=\sum_{i\le k\le p}{p-i\choose p-k}(x)_k/k!,
$$
we have 
$$
a_k=\sum_{i\le k\le p}w_i{p-i\choose p-k}/k!,
$$
implying that identity (\ref{eq3-8}) holds
and thus (\ref{sec7-e7}) follows.
\proofend

Theorem~\ref{w-sigmabar} is equivalent to the following result.
 
\begin{theo}
	\relabel{w-sigmabar2}
For any graph $G$ of order $p$, if 
\equ{w-sigmabar2-e1}
{\chi(G,x)=\sum_{i\le p}
a_i\, (x)_i,
}
then 
$$
w(G,x)=\sum_{i\le p} a_i i! x^i (1-x)^{p-i}.
$$
\end{theo}

\begin{cor}\relabel{cor3-2}
For any graph $G$ of order $p$, 
$$
w_k(G)=\sum_{i\le k}(-1)^{k-i} {p-i\choose p-k}i!a_i(G).
$$
where $a_i$ is given in 
(\ref{w-sigmabar2-e1}).
\end{cor}

\proof By Theorem~\ref{w-sigmabar}, 
\begin{align*}
w(G,x)&=(1-x)^p \bar \sigma(G,x/(1-x))
=(1-x)^p\sum_{i\le p}i!a_i (x/(1-x))^i
= \sum_{i\le p}i!a_i x^i(1-x)^{p-i}.
\end{align*}
Thus
$$
w_k=\sum_{i\le p}i! a_i{p-i\choose k-i}(-1)^{k-i}
=\sum_{i\le k}(-1)^{k-i} {p-i\choose p-k}i!a_i.
$$
\proofend

\begin{cor}\relabel{cor3-3} 
$$
w_k(N_p)=\sum_{i\le k}(-1)^{k-i} {p-i\choose p-k}i!S(p,i)
$$
and
$$
w(N_p,x)=
\sum _{i\le p} i!S(p,i)x^i(1-x)^{p-i}.
$$
\end{cor}

\proof As 
$$
a_i(N_p)=S(p,i),
$$
By Corollary~\ref{cor3-2},
$$
w_k(N_p)=\sum_{i\le k}(-1)^{k-i} {p-i\choose p-k}i!S(p,i).
$$
By Theorem~\ref{w-sigmabar2},
$$
w(N_p,x)=\sum_{i\le p}i!S(P,i)x^k(1-x)^{p-i}.
$$
\proofend

\subsection{Open problems}

Recall that a graph $G$ is said to be $w$-real if 
$w(G,x)$ has real roots only.
If $G$ is not $w$-real, it is called a $w$-unreal graph.

\begin{exa}
	 $C_4$ is the $w$-unreal graph with the minimal order:
	$$
	w(C_4,x)=2x^2(7x^2+4x+1).
	$$
\end{exa}

\begin{theo}[Brenti \cite{bre1992}]\relabel{th3-3}
$G$ is $w$-real 
if one of the following conditions is satisfied:
\begin{enumerate}[(b)]
\item $G$ contains a simplicial vertex $u$ 
such that $G-u$ is $w$-real;

\item (a special case of (a)) $G$ is chordal; and 

\item $G$ is the disjoin union of $G_1$ and $G_2$,
where each $G_i$ is $w$-real. 
\end{enumerate}
\end{theo}

\begin{que}\relabel{que3-3}
Find $w$-unreal graphs of order $5$.
\end{que}

\begin{con}[Brenti, Royle and Wagner \cite{bre1994}]\relabel{con3-3}
If both $G$ and $H$ are $w$-real and $G\cap H$ is complete, 
then $G\cup H$ is $w$-real.
\end{con}

\section{$\tau$-polynomial $\tau(G,x)$}
\relabel{sect-tau}

\subsection{Definition} 

Let $G$ be a simple graph of order $p$. 
For any $0\le i\le p$, let $c_i(G)$
(or simply $c_i$) denote the number in the following expression: 
\equ{sec8-e1}
{
\chi(G,x)=\sum_{0\le i\le p}(-1)^{p-i}c_i\langle x\rangle_i,
}
where $\angx_i=x(x+1)\cdots(x+i-1)$.

Now we define functions 
$\tau(G,x)$ and $\bar \tau(G,x)$   as follows
(see \cite{bre1992}):
\equ{sec8-e2}
{
\tau(G,x)=\sum_{0\le i\le p}c_i x^i
\quad \mbox{and}\quad 
\bar \tau(G,x)=\sum_{0\le i\le p}i! c_i x^i.
}

\begin{exa}
	\label{sec8-exa1}
\begin{enumerate}[(a)]
\item $\tau(K_1,x)=x$, as $\chi(K_1,x)=x=\angx_1$.

\item $\tau(K_2,x)=x^2+2x$, as 
$$
\chi(K_2,x)=x(x-1)=x(x+1-2)
=x(x+1)-2x=\angx_2-2\angx_1.
$$

\item $\tau(N_p,x)=B_p(x)$, 
where $
B_p(x)=\sum\limits_{1\le i\le p}S(p,i)x^i
$
is called a {\it Bell polynomial}.
\end{enumerate} 
\end{exa} 

\proof It is known that 
$$
\chi(N_p,x)=x^p=
\sum_{1\le i\le p}S(p,i)(x)_i.
$$
Letting $x=-z$ gives that 
$$
(-z)^p=\sum_{i\le p}S(p,i)(-z)_i=\sum_{i\le p}(-1)^iS(p,i) \langle z\rangle_i.
$$
Thus
$$
\chi(N_p,x)=x^p=\sum\limits_{1\le i\le p} (-1)^{p-i}S(p,i)\angx_i.
$$
By definition of $\tau(G,x)$, 
the conclusion follows.
\proofend

\subsection{Computation} 

If $G$ contains a simplicial 
vertex $u$, 
$\tau(G,x)$ can be expressed 
in terms of $\tau(G-u,x)$ 
and $\tau'(G-u,x)$.

\begin{pro}\relabel{sec8-prop1}
For any a simplicial vertex $u$ 
of $G$ with degree $k$,
\equ{sec8-e4}
{
\tau(G,x)=x\tau'(G-u,x)+(x+k)\tau(G-u,x).
}
\end{pro}

\proof By Theorem~\ref{sec5-th2}, 
$$
\sum_{i=0}^{p}(-1)^{p-i}
c_i\angx_i
=
\chi(G,x)=(x-k)\chi(G-u,x).
$$
Assume that 
$$
\chi(G-u,x)=\sum_{i=0}^{p-1}(-1)^{p-1-i}b_i\angx_i.
$$
Then 
\begin{eqnarray*}
(x-k)\chi(G-u,x)
&=&\sum_{i=0}^{p-1}(-1)^{p-1-i}b_i(x-k)\angx_i
\\&=&\sum_{i=0}^{p-1}(-1)^{p-1-i}b_i((x+i)-k-i)\angx_i
\\&=&\sum_{i=0}^{p-1}(-1)^{p-1-i}b_i\angx_{i+1}
+\sum_{i=0}^{p-1}(-1)^{p-i}b_i(k+i)\angx_i
\\&=&\sum_{j=1}^{p}(-1)^{p-j}b_{j-1}\angx_{j}
+\sum_{i=0}^{p-1}(-1)^{p-i}b_i(k+i)\angx_i.
\end{eqnarray*}
Thus, $c_i=b_{i-1}+b_i(k+i)$ for $i=0,1,\cdots,p-1,p$, 
where $b_p=0$. Hence 
\begin{eqnarray*}
\tau(G,x)=\sum_{i=0}^p c_ix^i 
&=&\sum_{i=0}^p (b_{i-1}+b_i(k+i))x^i
\\&=& x\tau(G-u,x)+k\sum_{i=0}^{p-1} b_ix^i
+x\sum_{i=1}^{p-1} ib_ix^{i-1}
\\&=&
x\tau(G-u,x)+k\tau(G-u,x)+x(\tau(G-u,x))'
\\&=&
(x+k)\tau(G-u,x)+x(\tau(G-u,x))'.
\end{eqnarray*}
\proofend


Now we are going to introduce 
the deletion-contraction 
formula for $\tau(G,x)$.

\begin{theo} \relabel{co-add-cont} 
For any graph $G$ with $e\in E(G)$, 
\equ{sec8-e5}
{
c_i(G)=c_i(G\backslash e)+c_i(G/e),
}
and 
\equ{sec8-e6}
{
\tau(G,x)=\tau(G\backslash e,x)+\tau(G/e,x),
}
where $G\backslash e$ (or resp. $G/e$) is  the 
graphs obtained from $G$ by removing $e$ 
(or resp. contracting $e$ and removing 
parallel edges but one). 
\end{theo} 

\proof (\ref{sec8-e5}) follows from the following identity:
\begin{eqnarray*}
\sum_{i\le p}(-1)^{p-i}
c_i\angx_i
&=&\chi(G,x)\\
&=&\chi(G\backslash e)
-\chi(G/e,x)\\ 
&=& \sum_{i\le p}(-1)^{p-i}c_i(G\backslash e)\angx_i
-\sum_{i\le p-1}(-1)^{p-1-i}c_i(G/e)\angx_i
\\ &=& \sum_{i\le p}(-1)^{p-i}[c_i(G\backslash e)+c_i(G/e)]\angx_i.
\end{eqnarray*}
By  (\ref{sec8-e5}),
 (\ref{sec8-e6}) is obtained.
\proofend

Thus, by Example~\ref{sec8-exa1} (c)
and Theorem~\ref{co-add-cont}, 
for any graph $G$,
$\tau(G,x)$ can be determined 
by applying
the following rules:
\begin{equation}
	\label{sec8-e7}
	\tau(G,x)=
	\left \{
	\begin{array}{ll}
		B_p(x), &\mbox{if }
		G\cong N_p;\\
\tau(G\backslash e,x)+\tau(G/e,x),
		\quad &\mbox{if }
		e\in E(G).
	\end{array}
	\right.
\end{equation}


\begin{exa} \label{sec8-exa2}
Let $P_2\cup K_1$ be the 
disjoint union of $P_2$ and $K_1$. 
Then applying
(\ref{sec8-e7})
yields that 
\begin{eqnarray*}
\tau(P_3,x)
&=&\tau(P_2\cup K_1,x)+\tau(P_2,x)
\\ &=& \tau(N_3,x)+2\tau(N_2,x)+\tau(N_1,x)
\\ &=& B_3(x)+2B_2(x)+B_1(x).
\end{eqnarray*}
\end{exa}

\begin{exa} \label{sec8-exa3}
Let $T$ be a tree of order $p$. Then
\equ{sec8-e8}
{
\tau(T,x)=\sum_{1\le k\le p}{p-1\choose k-1}B_{k}(x).
}
\end{exa}

\proof  The result is obvious 
when $p\le 2$.
Now assume that $p\ge 3$.
Let $u$ be a vertex in $T$
of degree $1$, and let $e$ be the edge incident with $u$. 
By (\ref{sec8-e7}), we have 
\equ{sec8-e9}
{
\tau(T,x)=
\tau(T-e,x)+\tau(T-u,x)
=\tau(K_1 \cup (T-u),x)+\tau(T-u,x).
}
As $T-u$ is a tree of order $p-1$, by the inductive assumption, we have 
\equ{sec8-e10}
{
	\tau(T-u,x)=\sum_{1\le k\le p-1}{p-2\choose k-1}B_{k}(x).
}
Thus, by (\ref{sec8-e9})
and (\ref{sec8-e10}), 
we have 
\Eqnn
{
\tau(T,x)
&=&\tau(K_1 \cup (T-u),x)+\tau(T-u,x)\\
&=&\sum_{1\le k\le p-1}{p-2\choose k-1}B_{k+1}(x)
+\sum_{1\le k\le p-1}{p-2\choose k-1}B_{k}(x)\\
&=&\sum_{2\le j\le p}
{p-2\choose j-2}
B_{j}(x)
+\sum_{1\le k\le p}
{p-2\choose k-1}B_{k}(x)\\
&=&\sum_{1\le k\le p}
{p-1\choose k-1}B_{k}(x).
}
\proofend 

\subsection{Properties}

\begin{pro}\relabel{pro4-4}
For any simple graph $G$ of order $p$,
if 
$$
\chi(G,x)=\sum_{k\le p} (-1)^{p-k}b_kx^k,
$$
then 
$$
\tau(G,x)=\sum_{k\le p}b_k B_k(x).
$$ 
\end{pro}

\proof Note that by Theorem~\ref{sec5-th4},  $b_k=\sum_{j}(-1)^{p-k+j}N_{j,k}$,
where $N_{j,k}$ is the number of spanning subgraphs 
of $G$ which have exactly $j$ edges 
and $k$ components.

It holds when $G=N_p$.
Then it can be proved by induction 
and applying 
Theorem~\ref{co-add-cont}.
Another proof is shown below
by the definition of  $\tau(G,x)$.
Note that 
\begin{eqnarray*}
\chi(G,x)&=&\sum_{k\le p}(-1)^{p-k}b_kx^k
\\ &=& \sum_{k\le p}(-1)^{p-k}b_k \sum_{i\le k}S(k,i)(-1)^{k-i}\angx_i
\\ &=& \sum_{i\le p}
(-1)^{p-i}
\left (\sum_{i\le k\le p}b_kS(k,i)\right )
\angx_i .
\end{eqnarray*}
Thus, by the definition of $\tau(G,x)$,
$$
\tau(G,x)=\sum_{i\le p}
\left (\sum_{i\le k\le p}b_kS(k,i)\right )
x^i
=\sum_{k\le p}b_k\sum_{i\le k}S(k,i)x^i
=\sum_{k\le p}b_kB_k(x).
$$
\proofend

\begin{cor}\relabel{cor4-0}
$$
\tau(K_p,x)=\sum_{k\le p}
\stirone{p}{k} 
B_k(x)
=\sum_{i\le p}\sum_{i\le k\le p}S(k,i)
\stirone{p}{k}x^i,
$$
where $\stirone{p}{k}$ 
is the {\it Stirling number of the first kind}, 
counting the number of permutations of $p$ elements 
with $k$ disjoint cycles.
\end{cor}

\subsection{Interpretations of coefficients}

Let $\Pi(G)$ be the set of partitions of $V(G)$.
For any $\P=\{P_1,\cdots,P_k\}
\in \Pi(G)$, let $G(\P)$ be the spanning subgraph of $G$ 
with edge set 
$\bigcup\limits_{1\le j\le k}
\{uv\in E(G): u,v\in P_j\}$.
Thus $G(\P)$ is the spanning subgraph obtained from $G$ by removing 
all edges whose ends are not in the same set $P_i$ of $\P$.

\begin{lem}[Brenti \cite{bre1992}]\relabel{lem4-1}
Let $G$ be a simple graph of order $p$.
For any $i\le p$, 
$$
c_i(G)=\sum_{\P\in \Pi(G) \atop |\P|=i} |\seta(G(\P))|,
$$
where $\seta(H)$ is the set of acyclic orientations of a graph $H$.
\end{lem}

\proof The result follows from
Theorem~\ref{co-add-cont}
and the fact that 
$$
|\seta(H)|=|\seta(H\backslash e)|+|\seta(H/e)|
$$
holds for any graph $H$ and edge $e$ in $H$ which is not a loop.
\proofend

Applying Lemma~\ref{lem4-1}, we have the following conclusions for $c_i(G)$:
For any simple graph $G$ of order $p$ and size $q$,
$$
c_0=0, c_1=|\seta(G)|, 
c_{p-1}={p\choose 2}+q,
c_p=1,
c_i\ge 1,
\forall i=2,3,\cdots,p-2
$$
and 
$$
c_{p-2}(G)={p\choose 3}+3{p\choose 4}+q{p-2\choose 2}
+m_2+\sum_{i=1}^3 l_i(2i-1),
$$
where $m_2$ is the number of matchings of $G$ with two edges
and $l_i$ is the number of induced subgraphs of $G$ with 
$3$ vertices and $i$ edges. 

\begin{theo}[Brenti \cite{bre1992}]\relabel{theo4-1}
For any graph $G$, 
$$
\tau(G,x)=\sum_{\P\in \Pi(G)}|\seta(G(\P))|x^{|\P|}.
$$
\end{theo}

\begin{cor}\relabel{cor4-1}
$$
\chi(G,x)=\sum_{\P\in \Pi(G)}|\seta(G(\P))|(-1)^{p-|\P|}\angx_{|\P|}.
$$
\end{cor}

\subsection{Open problems}

Recall that a graph $G$ is said to be $\tau$-real (resp. $\tau$-unreal)
if $\tau(G,x)$ has real roots only
(resp. has some roots not real).
Brenti, Royle and Wagner \cite{bre1994} obtained the following result.

\begin{theo}[Brenti, Royle and Wagner \cite{bre1994}]\relabel{th4-3}
$G$ is $\tau$-real  
if one of the following conditions is satisfied:
\begin{enumerate}[(a)]
\item $G$ is chordal;
\item $G$ is a cycle $C_p$, $p\ge 3$;
\item $G=H+ K_m$, where $H$ is $\tau$-real; 
\item $G=H+ K_m$, where $m$ is sufficiently large; and 
\item $G=G_1\cup G_2$, where each $G_i$ is $\tau$-real 
and $|V(G_1)\cap V(G_2)|\le 1$.
\end{enumerate}
\end{theo}

\begin{con}[Brenti, Royle and Wagner \cite{bre1994}]\relabel{con4-3}
If both $G$ and $H$ are $\tau$-real and $G\cap H$ is complete, 
then $G\cup H$ is $\tau$-real.
\end{con}

\begin{con}[Brenti, Royle and Wagner \cite{bre1994}]\relabel{con4-2}
Let $G$ and $H$ be vertex-disjoint graphs. 
If both $G$ and $H$ are $\tau$-real,
then the join $G + H$ is also $\tau$-real.
\end{con}

\begin{prob}[Brenti \cite{bre1992}]\relabel{prob4-1}
Is every graph $G$ $\tau$-real?
\end{prob}

So far no $\tau$-unreal graphs are known.

\def \ajc  
{{\it Australasian J. Combin.}}

\def \jcta {{\it J. Combin. Theory Ser. A}}
\def \jctb {{\it J. Combin. Theory Ser. B}}

\def \jgt {{\it J. Graph Theory}}

\def \cpc {{\it Combin. Probab. Comput.}}

\def \dm {{\it Discrete Math.}}

\def \eujc {{\it European J. Combin.}}

\def \ejc {{\it Electron. J. Combin.}}

\def \anc {{\it Ann. Combin.}}

\def \aam {{\it Advance in Applied Math.}}


\end{document}